\documentclass[leqno]{article}
\usepackage{latexsym,amsmath,amsfonts,mathrsfs,amssymb,amsthm}
\usepackage{yhmath}
\usepackage[usenames]{color}
\usepackage{graphicx}
\def\R{\mathbb R}
\def\N{\mathbb N}
\def\C{\mathbb C}
\def\Z{\mathbb Z}
\def\s{\sharp}
\def\a{\alpha}
\def\e{\epsilon}
\def\d{\delta}
\def\Y{\mathbb Y}
\def\T{\mathbb T}
\def\P{\mathbb P}
\def\be{\begin{equation}}
\def\ee{\end{equation}}
\def\wg{\wedge}
\def\bs{\backslash}
\def\qed{\hfill$\Box$\bigskip}
\def\nd{\noindent Proof. }
\def\H{{\cal H} }
\numberwithin{equation}{section}
\newtheorem{lem}[equation]{Lemma}
\newtheorem{pro}[equation]{Proposition}
\newtheorem{defn}[equation]{Definition}
\newtheorem{thm}[equation]{Theorem}

\newtheorem{cor}[equation]{Corollary}

\newtheorem{rem}[equation]{Remark}
\begin{document}

\begin{center} \Large\textbf{On the topological minimality of unions of planes of arbitrary dimension}\end{center}

\bigskip

\centerline{\large Xiangyu Liang}

\vskip 1cm

\centerline {\large\textbf{Abstract.}}In this article we prove the topological minimality of unions of several almost orthogonal planes of arbitrary dimensions. A particular case was proved in \cite{2p}, where we proved the Almgren minimality (which is a weaker property than the topological minimality) of the union of two almost orthogonal 2 dimensional planes. On the one hand, the topological minimality is always proved by variations of calibration methods, but in this article, we give a continuous family topological minimal sets, hence calibrations cannot apply. The advantage of a set being topological minimal (compared to Almgren minimal) is that its product with $\R^n$ stays topological minimal. This leads also to finding minimal sets which are unions of non transversal (hence far from almost orthogonal) planes; On the other hand, regularity for higher dimensional minimal sets is much less clear than those of dimension 2, hence more efforts are needed for higher dimensional cases. 

\bigskip


\bigskip


\bigskip


\section{Introduction and preliminaries}

\subsection{Introduction}

In this article we discuss topological minimality of unions of planes of arbitrary dimensions. A particular case was proved in \cite{2p}, where we proved the Almgren minimality (which is a weaker property than the topological minimality) of the union of two almost orthogonal 2 dimensional planes. 

The notion of Almgren-minimality (introduced in \cite{Al76}) is a general notion of weak solutions, in the setting of sets, of Plateau's problem, which aims at understanding the regularity and existence of physical objects that have certain minimizing properties such as soap films. Roughly speaking, we say that a closed set $E$ is $d$-dimensional Almgren-minimal when there is no deformation $F=\varphi(E)$, where $\varphi$ is Lipschitz and $\varphi(x)-x$ is compactly supported, for which the Hausdorff measure $\H^d(F)$ is smaller than $\H^d(E)$. See Definition 1.11 for the  precise definition.

The notion of topological minimal sets (introduced by the author in \cite{topo}) is also in the setting of sets, but instead of minimizing Hausdorff measure among compact deformations, one asks that a topological minimal set admits a minimal measure among all sets that satisfy some topological property. At first glance, topological minimality is stronger than Almgren minimality (See Proposition 1.18), though we do not know whether it is strictly stronger. The definition of topological minimal sets might be physically less intuitive than that of Almgren minimal sets, however one of its big advantages is that one always have existence for topological minimal sets (cf. \cite{topo} Theorems 4.2 and 4.28), while for Almgren minimal sets, one can only expect for partial results. Also, we know that the product of a topological minimal set with $\R^n,n\in\N$ is still topologically minimal (cf. \cite{topo} Proposition 3.23), but we do not have a such property for Almgren minimal sets.

\bigskip

Compared to chains, currents, rectifiable varifolds of sets of least perimeter, etc., which are more usually used to modernize Plateau's problem, minimal sets has less structure, and very little is known for their regularity. On the other hand, minimal sets are absolute minimizers, and thus we expects better regularity than for general critical points. 
%

\medskip

Now let us say something more about regularity of minimal sets. Since topological minimal sets are automatically Almgren minimal, all the regularity properties stated below also hold for topological minimizers.

First regularity results for Almgren-minimal sets have been given by Frederick Almgren \cite{Al76} (rectifiability, Ahlfors regularity in arbitrary dimension), then generalized by Guy David and Stephen Semmes \cite{DS00} (uniform rectifiability, big pieces of Lipschitz graphs). 

Since Almgren minimal sets are rectifiable and Ahlfors regular, they admit a tangent plane at almost every point. But our main interest is to study those points where there is no tangent plane, i.e. singular points. In \cite{DJT}, Guy David proved that at each point of an Almgren minimal set,  every blow-up limit (see Definition 1.23) is an Almgren minimal cone, that is, an Almgren minimal set which is a cone at the same time (we will call them minimal cones throughout the paper, since topological minimal cones will not be mentioned). Thus, the study of singular points is transformed into the classification of singularities, i.e., into looking for a list of minimal cones. 

In $\R^3$, the list of minimal cones has been given by several mathematicians a century ago. (See for example \cite{La} or \cite{He}). 
%
For example, 2-dimensional minimal cones in $\R^3$ are, modulo isomorphism: a plane, a $\Y$ set (the union of 3 half planes that meet along a straight line where they make angles of 120 degrees), and a $\T$ set (the cone over the 1-skeleton of a regular tetrahedron centered at the origin). See the pictures below.

\centerline{\includegraphics[width=0.2\textwidth]{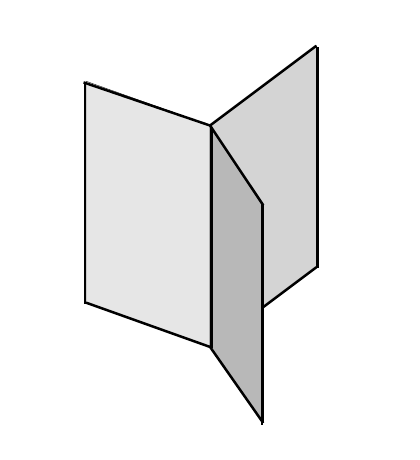} \hskip 2cm\includegraphics[width=0.25\textwidth]{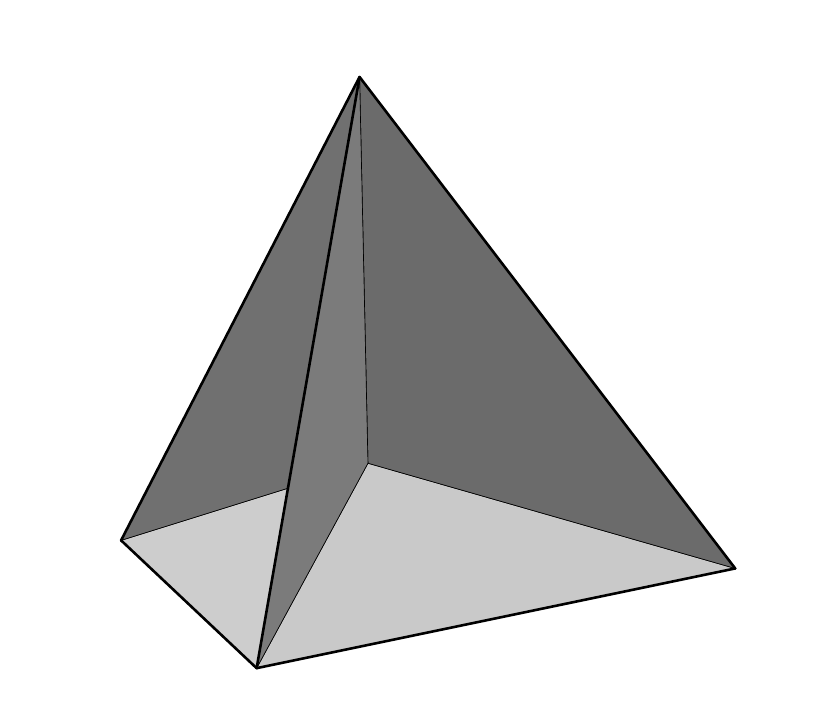}}

In higher dimensions, even in dimension 4, the list of minimal cones is still very far from clear. Up to now we only know some particular example, such as the cone over the $n-2$ skeleton of a regular simplex centered at the origin in $\R^n$ for $n\ge 2$ (\cite{LM94}), the cone over the $n-2$ dimensional skeleton of cubes centered at the origin in $\R^n$ for $n\ge 4$ (\cite{Br91}), the almost orthogonal union of two planes (\cite{2p}), the set $Y\times Y$ which is the product of two 1-dimensional $\Y$ sets (\cite{YXY}).

Among all the above minimal cones, the minimality of most of them are proved by calibrations (or some generalized calibrations). Essencially,  all sorts of calibration methods always prove directly the topological minimality, rather than the weaker Almgren minimality. In this case, from all those calibrated minimal sets, we can obtain new higher dimensional minimal cones by simply taking their products with $\R^n$.

However for the unions of two almost orthogonal planes, to the author's knowledge, no calibration works for them, and the proof is very different. (The non existence of calibration might also be the reason why we can have such a continuous family of minimal cones.) However it still makes sense to ask whether they are topologically minimal or not, for example this is related to the interesting question that whether there exists unions of non-transversal planes that are minimal. (This can not happen for 2-dimensional case, see for example \cite{DJT} Proposition 14.1). The easist way to get an affirmative answer is to prove that in fact the unions of almost orthogonal planes are also topologically minimal (even though they are not well calibrated), and then the product of $\R^k$ with theses unions will be minimal cones, which are unions of non-transversal planes.

On the other hand, since the proof of minimality of union of two almost orthogonal planes is very different, it is natural to ask whether we can also prove the result for more general cases,  i.e. the union of more than two almost orthogonal higher dimensional planes. 

Under the above two motivation, we will discuss the topological minimality of the union of several almost orthogonal planes of dimension $d\ge 2$ in this paper. We will prove the following :

\begin{thm}[Topological minimality of the union of $m$ almost orthogonal planes]\label{main}For each $d\ge 2$ and $m\ge 2$, there exists $\theta=\theta(m,d)\in]0,\frac\pi2[$, such that if $P^1,P^2,\cdots,P^m$ are $m$ planes of dimension $d$ in $\R^{dm}$ with characteristic angles $\alpha^{ij}=(\alpha^{ij}_1,\alpha^{ij}_2,\cdots,\a^{ij}_d)$ between $P^i$ and $P^j,1\le i<j\le m$, which verify $\theta<\a^{ij}_1\le\a^{ij}_2\le\cdots\le\a^{ij}_d\le\frac\pi2$ for all $1\le i<j\le m$, then their union $\cup_{i=1}^m P^i$ is a topological minimal cone.\end{thm}

The characteristic angles of two tansversal $d-$planes $P^1,P^2$ is $\alpha=(\alpha_1,\alpha_2,\cdots,\a_d)$ implies that there exists an orthonormal basis $\{e_i\}_{1\le i\le 2m}$ of $\R^{2m}$ (the linear subspace generated by $P^1$ and $P^2$) such that $P^1$ is generated by $\{e_i\}_{1\le i\le m}$, and $P^2$ is generated by $\{\cos\a_i e_i+\sin\a_i e_{n+i}\}$, see Definition 2.5 for a precise definition. Hence the characteristic angles describe their relative position. An almost equivalent statement of Theorem \ref{main} that might be easier to understand is that any union of  $m$ $d$-dimensional planes which are mutually almost orthogonal in $\R^{dm}$ is topologically minimal, i.e. there exists $a>0$ (small), such that if for any $1\le i<j\le m$, and any $v_1\in P^i,v_2\in P^j$,
\be |<v_1,v_2>|\le a||v_1||||v_2||\ee 
then $\cup_{1\le i\le m}P^i$ is minimal. Or intuitively, there is an``open set" of unions of $m$ planes, which contains the orthogonal union of $m$ planes, such that each element in this set (which is a union of $m$ planes) is a topological minimal cone.

Note that when the angles between planes are small, their union cannot be Almgren minimal, because we can easily "pinch" two planes in the center and decrease measure. See the construction in \cite{La94}.

As a corollary of Theorm 1.11, we will have families of unions of non transversal $d$-planes $d>2$ which are minimal, by simply taking the product of $\R^{d-k}$ with the union of $m$ almost orthogonal $k$-planes, with $2\le k< d$. See Corollary 9.2.

%

The general plan for the proof for Theorem \ref{main} will be similar to that in \cite{2p}, but due to the lack of knowledge of regularity for higher dimensional minimal sets, as well as the difference between Almgren minimal sets and topological minimal sets, there are substantial technical differences: in particular, the uniqueness theorem of the orthogonal union (Thm 3.1) of higher dimensional planes, and the projection property for topological minimal competitors (Proposition 6.1). The treat of harmonic extensions is also different in higher dimensions. The existence of minimal topological competitors (Theorem 4.5) is also different from the partial existence result used in \cite{2p}, but this was already proved in \cite{topo}.

The rest of this paper will be organized as the following. We will give details at places where proofs are different (Section 2,3,7, part of Section 6), and will only sketch the prove for the rest part.

In Subsection 1.1 we will give some basic definition and notation that we will use frequently afterwards.

Section 2 will be devoted to estimate the sum of projections of unit simple $d-$vectors on several $d$-planes, depending on their mutual characteristic angles.  Based on this we will estimate the sum of the measure of the projections of a rectifiable set.

Section 3 will be devoted to prove the uniqueness theorem.

In Section 4 and 5 we will sketch the construction of converging sequences $E_k$ of topologically minimal competitors for $P_k$ (where $P_k$ is a sequence of unions of planes that converges to the orthogonal union of planes), and the construction of the $\d$-process.

We prove some necessary regularity results for $E_k$ in Section 6.  

In Section 7 we give the estimate of Dirichlet energy of a graph on $d-$dimensional annulus, with prescribed boundary condition.

The proof of Thm \ref{main} will be given in Section 8.

Finally in Section 9, we will give an corollary about minimality of unions of non transversal higher dimensional planes, and discuss some related open problems.

\subsection{Preliminaries}

$[a,b]$ is the line segment with end points $a$ and $b$;

$\overrightarrow{ab}$ is the vector $b-a$; while being specified, it can also represent the half line issued from the point $a$ and passing through $b$;

$B(x,r)$ is the open ball with radius $r$ and centered on $x$;

$\overline B(x,r)$ is the closed ball with radius $r$ and center $x$;

$\H^d$ is the Hausdorff measure of dimension $d$ ;

$d_H(E,F)=\max\{\sup\{d(y,F):y\in E\},\sup\{d(y,E):y\in F\}\}$ is the Hausdorff distance between two sets $E$ and $F$.

$d_{x,r}$ : while not being specified, it denotes the relative distance with respect to the ball $B(x,r)$, defined by
$$ d_{x,r}(E,F)=\frac 1r\max\{\sup\{d(y,F):y\in E\cap B(x,r)\},\sup\{d(y,E):y\in F\cap B(x,r)\}\}.$$

In the next definitions, fix integers $0<d<n$. We first give a general definition for minimal sets. Briefly, a minimal set is a closed set which minimizes the Hausdorff measure among a certain class of competitors. Different choices of classes of competitors give different kinds of minimal sets.

\begin{defn}[Minimal sets]Let $0<d<n$ be integers. Let $U\subset \R^n$ be an open set. A relative closed set $E\subset U$ is said to be minimal of dimension $d$ in $U$ with respect to the competitor class $\mathscr F$ (which contains $E$) if 
\be \H^d(E\cap B)<\infty\mbox{ for every compact ball }B\subset U,\ee
and
\be \H^d(E\bs F)\le \H^d(F\bs E)\ee
for any competitor $F\in\mathscr F$.
\end{defn}

\begin{defn}[Almgren competitor (Al competitor for short)] Let $E$ be relatively closed in an open subset $U$ of $\R^n$. An Almgren competitor for $E$ is an relatively closed set $F\subset U$ that can be written as $F=\varphi_1(E)$, where $\varphi_t:U\to U,t\in [0,1]$ is a family of continuous mappings such that 
\be \varphi_0(x)=x\mbox{ for }x\in U;\ee
\be\mbox{ the mapping }(t,x)\to\varphi_t(x)\mbox{ of }[0,1]\times U\mbox{ to }U\mbox{ is continuous;}\ee
\be\varphi_1\mbox{ is Lipschitz,}\ee
  and if we set $W_t=\{x\in U\ ;\ \varphi_t(x)\ne x\}$ and $\widehat W=\bigcup_{t\in[0.1]}[W_t\cup\varphi_t(W_t)]$,
then
\be \widehat W\mbox{ is relatively compact in }U.\ee
 
Such a $\varphi_1$ is called a deformation in $U$, and $F$ is also called a deformation of $E$ in $U$.
\end{defn}

\begin{defn}[Almgren minimal sets]
Let $0<d<n$ be integers, $U$ be an open set of $\R^n$. An Almgren-minimal set $E$ in $U$ is a minimal set defined in Definition 1.3 while taking the competitor class $\mathscr F$ to be the class of all Almgren competitors for $E$.\end{defn}

\begin{rem}When the ambient set $U$ is $\R^n$, or a ball, we can also take the class of local Almgren competitors to define the same notion of minimal set. Keep the $E$, $U$, $n$ and $d$ as before; a local Almgren competitor of $E$ in $U$ is a set $F=f(E)$, with
\be f=id\mbox{ outside some compact ball }B\subset U,\ee
\be f(B)\subset B,\ee 
and $f$ is Lipschitz. 

A such $f$ is called a local deformation in $U$, or a deformation in $B$, and $F=f(E)$ is also called a local deformation of $E$ in $U$, or a deformation of $E$ in $B$.

Note that in this case, the condition (1.5) becomes 
\be \H^d(E\cap B)\le \H^d(F\cap B).\ee

We say that a set $E$ closed in an open set $U$ is locally minimal if (1.4) holds, and for any compact ball $B\subset U$, and any local Almgren competitor $F$ for $E$ in $B$, (1.15) holds.

One can easily verify that when $U$ is $\R^n$ or a ball, the class of Al competitors coincides with the class of local Al competitors, so the two classes define the same kind of minimal sets. However, if the ambient set $U$ has a more complicated geometry, then the class of local Al competitors is strictly smaller, so a set minimizing the Hausdorff measure among local Al competitors might fail to be Al-minimal. 
\end{rem}

\begin{rem}In general, the notion of minimal sets does not depend much on the ambient dimension. For example one can easily check that $E\subset U$ is $d-$dimensional Almgren minimal in $U\subset \R^n$ if and only if $E$ is Almgren minimal in $U\times\R^m\subset\R^{m+n}$, for any integer $m$.\end{rem}

\begin{defn}[Topological competitors] Let $E$ be a closed set in $\R^n$. We say that a closed set $F$ is a topological competitor of dimension $d$ ($d<n$) of $E$, if there exists a ball $B\subset\R^n$ such that

1) $F\bs B=E\bs B$;

2) For all Euclidean $n-d-1$-sphere $S\subset\R^n\bs(B\cup E)$, if $S$ represents a non-zero element in the singular homology group $H_{n-d-1}(\R^n\bs E;\Z)$, then it is also non-zero in $H_{n-d-1}(\R^n\bs F;\Z)$.
We also say that $F$ is a topological competitor of $E$ in $B$.
\end{defn}

And Definition 1.3 gives the definition of topological minimizers.

The simplest example of a topological minimal set is a $d-$dimensional plane in $\R^n$.  

\begin{pro}[cf.\cite{topo} Corollary 3.17] All topological minimal sets are Almgren minimal.
\end{pro}

\begin{defn}[reduced set] Let $U\subset \R^n$ be an open set. For every closed subset $E$ of $U$, denote by
\be E^*=\{x\in E\ ;\ \H^d(E\cap B(x,r))>0\mbox{ for all }r>0\}\ee
 the closed support (in $U$) of the restriction of $\H^d$ to $E$. We say that $E$ is reduced if $E=E^*$.
\end{defn}

It is easy to see that
\be \H^d(E\bs E^*)=0.\ee
In fact we can cover $E\bs E^*$ by countably many balls $B_j$ such that $\H^d(E\cap B_j)=0$.

\begin{rem}It is not hard to see that if $E$ is Almgren minimal (resp. topological minimal), then $E^*$ is also Almgren minimal (resp. topological minimal). As a result it is enough to study reduced minimal sets.
\end{rem}

\begin{defn}[blow-up limit] Let $U\subset\R^n$ be an open set, let $E$ be a relatively closed set in $U$, and let $x\in E$. Denote by $E(r,x)=r^{-1}(E-x)$. A set $C$ is said to be a blow-up limit of $E$ at $x$ if there exists a sequence of numbers $r_n$, with $\lim_{n\to \infty}r_n=0$, such that the sequence of sets $E(r_n,x)$ converges to $C$ for the Hausdorff distance in any compact set of $\R^n$.
\end{defn}

\begin{rem}A set $E$ might have more than one blow-up limit at a point $x$. However it is not known yet whether this can happen to minimal sets.
\end{rem}

\begin{pro}[c.f. \cite{DJT} Proposition 7.31]Let $E$ be a reduced Almgren minimal set in the open set $U$, and let $x\in E$. Then every blow-up limit of $E$ at $x$ is a reduced Almgren minimal cone $F$ centred at the origin, and $\H^d(F\cap B(0,1)=\theta(x):=\lim_{r\to 0} r^{-d}\H^d(E\cap B(x,r)).$\end{pro}

An Almgren minimal cone is just a cone which is also Almgren minimal. We will call them minimal cones throughout this paper, since we will not talk about any other type of minimal cones. Also, when not specified, minimal set will mean Almgren minimal set in the rest of the paper.

\begin{rem}The existence of the density $\theta(x)$ is due to the monotonicity of the density function $\theta(x,r):=r^{-d}\H^d(E\cap B(x,r))$ for minimal sets. See for example \cite{DJT} Proposition 5.16.
\end{rem}

We now state some regularity results on 2-dimensional Almgren minimal sets. Note that these properties also holds for any topological minimal set, after Proposition 1.18.

\begin{defn}[bi-H\"older ball for closed sets] Let $E$ be a closed set of Hausdorff dimension 2 in $\R^n$. We say that $B(0,1)$ is a bi-H\"older ball for $E$, with constant $\tau\in(0,1)$, if we can find a 2-dimensional minimal cone $Z$ in $\R^n$ centered at 0, and $f:B(0,2)\to\R^n$ with the following properties:

$1^\circ$ $f(0)=0$ and $|f(x)-x|\le\tau$ for $x\in B(0,2);$

$2^\circ$ $(1-\tau)|x-y|^{1+\tau}\le|f(x)-f(y)|\le(1+\tau)|x-y|^{1-\tau}$ for $x,y\in B(0,2)$;

$3^\circ$ $B(0,2-\tau)\subset f(B(0,2))$;

$4^\circ$ $E\cap B(0,2-\tau)\subset f(Z\cap B(0,2))\subset E.$  

We also say that B(0,1) is of type $Z$.

We say that $B(x,r)$ is a bi-H\"older ball for $E$ of type $Z$ (with the same parameters) when $B(0,1)$ is a bi-H\"older ball of type $Z$ for $r^{-1}(E-x)$.
\end{defn}

\begin{thm}[Bi-H\"older regularity for 2-dimensional Almgren minimal sets, c.f.\cite{DJT} Thm 16.1]\label{holder} Let $U$ be an open set in $\R^n$ and $E$ a reduced Almgren minimal set in $U$. Then for each $x_0\in E$ and every choice of $\tau\in(0,1)$, there is an $r_0>0$ and a minimal cone $Z$ such that $B(x_0,r_0)$ is a bi-H\"older ball of type $Z$ for $E$, with constant $\tau$. Moreover, $Z$ is a blow-up limit of $E$ at $x$.
\end{thm}

\begin{defn}[point of type $Z$] In the above theorem, we say that $x_0$ is a point of type $Z$ (or $Z$ point for short) of the minimal set $E$.
\end{defn}

\begin{rem} Again, since we might have more than one blow-up limit for a minimal set $E$ at a point $x_0\in E$, the point $x_0$ might of more than one type (but all blow-up limits at a point are bi-H\"older equivalent). However, if one of the blow-up limits of $E$ at $x_0$ admits the``full-length'' property (see Remark \ref{ful}), then in fact $E$ admits a unique blow-up limit at the point $x_0$. Moreover, we have the following $C^{1,\a}$ regularity around the point $x_0$. In particular, the blow-up limit of $E$ at $x_0$ is in fact a tangent cone of $E$ at $x_0$.
\end{rem}

\begin{thm}[$C^{1,\a}-$regularity for 2-dimensional minimal sets, c.f. \cite{DEpi} Thm 1.15]\label{c1} Let $E$ be a 2-dimensional reduced minimal set in the open set $U\subset\R^n$. Let $x\in E$ be given. Suppose in addition that some blow-up limit of $E$ at $x$ is a full length minimal cone (see Remark \ref{ful}). Then there is a unique blow-up limit $X$ of $E$ at $x$, and $x+X$ is tangent to $E$ at $x$. In addition, there is a radius $r_0>0$ such that, for $0<r<r_0$, there is a $C^{1,\a}$ diffeomorphism (for some $\a>0$) $\Phi:B(0,2r)\to\Phi(B(0,2r))$, such that $\Phi(0)=x$ and $|\Phi(y)-x-y|\le 10^{-2}r$ for $y\in B(0,2r)$, and $E\cap B(x,r)=\Phi(X)\cap B(x,r).$ 

We can also ask that $D\Phi(0)=Id$. We call $B(x,r)$ a $C^1$ ball for $E$ of type $X$.
\end{thm}

\begin{rem}[full length, union of two full length cones $X_1\cup X_2$]\label{ful}We are not going to give the precise definition of the full length property. Instead, we just give some information here, which is enough for the proofs in this paper.

1) The three types of 2-dimensional minimal cones in $\R^3$, i.e. the planes, the $\Y$ sets, and the $\T$ sets, all verify the full-length property (c.f., \cite{DEpi} Lemmas 14.4, 14.6 and 14.27). Hence all 2-dimensional minimal sets $E$ in an open set $U\subset\R^3$ admits the local $C^{1,\a}$ regularity at every point $x\in E$. But this was known from \cite{Ta}.

2) (c.f., \cite{DEpi} Remark 14.40) Let $n>3$. Note that the planes, the $\Y$ sets and the $\T$ sets are also minimal cones in $\R^n$. Denote by $\mathfrak C$ the set of all planes, $\Y$ sets and $\T$ sets in $\R^n$. Let $X=\cup_{1\le i\le n}X_i\in \R^n$ be a minimal cone, where $X_i\in \mathfrak{C}, 1\le i\le n$, and for any $i\ne j$, $X_i\cap X_j=\{0\}$. Then $X$ also verifies the full-length property. 
\end{rem}

\begin{lem}[Structure of 2-dimensional minimal cones in $\R^n$, cf. \cite{DJT} Proposition 14.1] Let $X$ be a 2-dimensional minimal cone in $\R^n$, and set $K=X\cap\partial B(0,1)$. Then $K$ is a finite union of great circles or arcs of great circles $C_j, j\in J$. The $C_j$  can only meet at their extremities, and each extremity is a common extrimity of exactly three $C_j$, which meet with $120^\circ$. 
\end{lem}

As for the regularity for minimal sets of higher dimensions, we know much less. But for points which admit a tangent plane (i.e. some blow up-limit on the point is a plane), we still have the $C^1$ regularity.

\begin{pro}[cf.\cite{2p} Proposition 6.14]\label{e1}For $2\le d<n<\infty$, there exists $\epsilon_1=\e_1(n,d)>0$ such that if $E$ is a $d$-dimensional minimal set in an open set $U\subset\R^n$, with $B(0,2)\subset U$ and $0\in E$. Then if $E$ is $\epsilon_1$ near a $d-$plane $P$ in $B(0,1)$, then $E$ coincides with the graph of a $C^1$ map $f:P\to P^\perp$ in $B(0,\frac 34)$. Moreover $||\nabla f||_\infty<1$.
\end{pro}

\begin{rem}This proposition is a direct corollary of Allard's famous regularity theorem for stationary varifold. See \cite{All72}.
\end{rem}

%

\medskip

\section{Some basic preliminaries and estimates for unit simple $d$-vectors}
Let $2\le d<n$ be two integers. Denote by $\wg_d(\R^n)$ the space of all $d$-vectors in $\R^n$. 
Set $I_{n,d}=\{I=(i_1,i_2,\cdots, i_d):1\le i_1<i_2<\cdots<i_d\le n\}$.
Let $\{e_i\}_{1\le i\le n}$ be an orthonormal basis of $\R^n$. For any $I=(i_1,i_2,\cdots, i_d)\in I_{n,d}$, denote by $e_I=e_{i_1}\wg e_{i_2}\wg\cdots\wg e_{i_d}$. Then the set $\{e_I,I\in I_{n,d}\}$ forms a basis of $\wg_d(\R^n)$. The standard scalar product on $\wg_d(\R^n)$ is: for $\xi=\sum_{I\in I_{n,d}}a_I e_I$ and $\zeta=\sum_{I\in I_{n,d}}b_I e_I$,
\be <\xi,\zeta>=\sum_{I\in I_{n,d}}a_Ib_I.\ee

Denote by $|\cdot|$ the norm induced by this scalar product.

Now given a unit simple $d$-vector $\xi$, we can associate it to a $d$-dimensional subspace $P(\xi)\in G(n,d)$, where $G(n,d)$ denotes the set of all $d$-dimensional subspace of $\R^n$:
\be P(\xi)=\{v\in\R^n,v\wedge\xi=0.\}\ee
In other words, if $\xi=x_1\wg x_2\cdots\wg x_d$, $x_1,\cdots x_d$ being orthogonal, then $P(\xi)$ is the $d-$subspace generated by $\{x_i\}_{1\le i\le d}$.

From time to time, when there is no ambiguity, we also write $P=x_1\wg x_2\wg\cdots x_d$, where $P\in G(n,d)$ and $\{x_i\}_{1\le i\le d}$  are $d$ unit vectors such that $P=P(x_1\wg\cdots \wg x_d)$.

Now if $f$ is a linear map from $\R^n$ to $\R^n$, then we denote by $\wg_d f$ (and sometimes by $f$ if there is no ambiguity) the linear map from $\wg_d(\R^n)$ to $\wg_d(\R^n)$ such that
\be \wg_df(x_1\wg x_2\cdots\wg x_d)=f(x_1)\wg f(x_2)\wg\cdots\wg f(x_d).\ee

And if $P\in G(n,d)$, then $P=P(\xi)$ for some unit simple $d$-vector $\xi$ (such a $d$ vector always exists), we define $|f(\cdot)|:G(n,d)\to\R^+\cup\{0\}$ by
\be |f(P)|=|\wg_df(x_1\wg\cdots \wg x_d).\ee
One can easily verify that the value of $|f(P)|$ does not depend on the choice of the unit simple vector $\xi$ that generates $P$. Hence $|f(\cdot)|$ is well defined.

Now let us recall the definition of characteristic angles between two $d-$planes. 

\begin{defn} Let $P^1,P^2$ be two $d$-dimensional planes in $\R^n$. Among all pairs of unit vectors $(v,w)$ with $v\in P^1,\ w\in P^2$, we choose $(v_1,w_1)$ which minimizes the angle between them. We denote by $\a_1$ this angle. Next we look at all the pairs of unit vectors $\{(v',w'):v'\in P^1,w'\in P^2,v'\perp v_1,w'\perp w_1\}$ , and we choose $(v_2,w_2)$ which minimizes the angle among all such pairs. Denote by $\a_2$ this angle. We continue like this, and then we get $d$ angles $\a_1\le\a_2\le\cdots\a_d$. They are the $d$ characteristic angles of $P^1$ and $P^2$. Or alternatively, we call the $d$-tuple $\a=(\a_1,\cdots \a_d)$ (with $\a_1\le\a_2\le\cdots\a_d$) the characteristic angle between $P^1$ and $P^2$, and $\min\a=\a_1$ the smallest angle in this $d-$tuple.
\end{defn}

Characteristic angles characterize absolutely the relative position between two planes, in the sense that we can find an orthonormal basis $\{e_i\}_{1\le i\le n}$ of $\R^n$, such that
\be P^1=P(e_1\wedge e_2\cdots\wg e_d)\mbox{ and }P^2=P[\bigwedge_{i=1}^d(\cos\a_i e_i+\sin\a_i e_{d+i})].\ee

Two $d-$planes are said to be orthogonal to each other if their characteristic angles are all $\frac\pi 2$.

We are now going to give some estimates on projections of an $d-$vector on $d-$planes. First comes the orthogonal case.

 \begin{lem}\label{F1}(c.f. \cite{Mo84} Lemma 5.2)

Let $P,Q$ be two subspaces of $\R^n$ with 
\be dim(P\cap {Q}^\perp)\ge dim\ P-d+2\ee 
Let $\xi$ be a simple unit $d-$vector in 
$\wg_d(\R^n)$. Denote by $p,q$ the orthogonal projections from $\R^n$ onto $P$ and $Q$ respectively. Then the projections of $\xi$ verify
\be|p(\xi)|+|q(\xi)|\le 1.\ee
Moreover, if 
\be dim(P\cap Q)<d-2,\ee 
then
\be |p\xi|+|q\xi|=1\mbox{ if and only if }\xi\mbox{ belongs to }P\mbox{ or }Q.\ee
\end{lem}

As a corollary of this, we have

\begin{pro}Let $d\ge 2$, and $E_1,E_2$ two Almgren minimal sets of dimension $d$ in $\R^{m_1}$ and $\R^{m_2}$ respectively. Then the orthogonal union  $E_1\cup E_2$ is an Almgren minimal set in 
$\R^{m_1+m_2}$.
\end{pro}

\nd Let $F$ be a deformation of $E_1\cup E_2$ in $\R^{m_1+m_2}$, then there exists $R>0$ and $f$ a Lipschitz deformation in $\R^{m_1+m_2}$ such that  
\be f(B(0,R))\subset B(0,R);f|_{B(0,R)^C}=Id,\mbox{ and }f(E_1\cup E_2)=F.\ee

Denote by $p_i$ the projection on $\R^{m_i}$, $i=1,2$. Then $p_i\circ f(E_i)$ is a deformation of $E_i$ in $B(0,R)\cap \R^{m_i}$, $i=1,2$. By the Almgren minimality of $E_i$, $\H^d(p_i\circ f(E_i))\ge H^d(E_i)$, hence
\be \H^d(p_i(E))=\H^d(p_i\circ f(E_1\cup E_2))\ge \H^d(p_i\circ f(E_i))\ge \H^d (E_i),i=1,2.\ee

Then we apply Lemma \ref{F1}, and the following lemma, we obtain that
\be \H^d(E)\ge \H^d(p_1(E))+\H^d(p_2(E))\ge \H^d(E_1)+\H^d(E_2)=\H^d(E_1\cup E_2),\ee
where the conclusion follows.\qed

\begin{lem}
 Let $n>d\ge 2$, and $P,Q$ be two subspaces in $\R^n$, $F\subset \R^n$ be a $d-$rectifiable set. Denote by $p,q$ the orthogonal projections on $P$ and $Q$ respectively. Let $\lambda\ge 0$ be such that for almost all $x\in F$, the tangent plane $T_xF\in G(n,d)$ of $F$ verifies
\be |p(T_xF)|+|q(T_xF)|\le\lambda.\ee
Then
\be \H^d(p(F))+\H^d(q(F))\le\lambda H^d(F).\ee
\end{lem}

\nd

Denote by $f$ the restriction of $p$ on $F$, then $f$ is a Lipschitz function from a $d$-rectifiable set to a $d-$rectifiable subset of $P$. Since $F$ is $d$-rectifiable, for $\H^d-$almost all $x\in F$, $f$ has an approximate differential
\be ap Df(x): T_xF\to P \ee (c.f.\cite{Fe}, Thm 3.2.19). Moreover this differential is such that $||\bigwedge_d ap Df(x)||\le 1$ almost everywhere, because $f$ is $1-$Lipschitz.

Now we can apply the area formula to $f$, (c.f. \cite{Fe} Cor 3.2.20). For all $\H^d|_F$-integrable functions $g\ : \ F\to\bar\R$, we have
\be \int_F (g\circ f)\cdot ||\wedge_d apDf(x)||d\H^d=\int_{f(F)} g(z)N(f,z)d\H^dz,\ee where $ N(f,z)=\s\{f^{-1}(z)\},$
and for $z\in p(F)$ we have $N(f,z)\ge 1$.
Take $g\equiv1$, we get
\be\int_F ||\wedge_d apDf(x)||d\H^d=\int_{p(F)} N(f,z)d\H^dz\ge \int_{p(F)}d\H^d=\H^d(p(F)).\ee

Recall that $p$ is linear, hence its differential is itself. As a result $apDf(x)$ is the restriction of $p$ on the $d$-subspace $T_xF$, which implies that if $\{u,v\}$ is an orthonormal basis of $T_xF$, then
\be ||\wedge_dapDf(x)||=|p(T_xF)|\ee
by (2.4).
 Hence by (2.21)
\be\int_F |p(T_xF)| d\H^d(x)\ge \H^d(p(F)).\ee
A similar argument gives also:
\be\int_F |q(T_xF)| d\H^d(x)\ge \H^d(q(F)).\ee
Summing (2.23) and (2.24) we get
\be\begin{array}{ll} \H^d(p F)+\H^d(q F)&\le\int_F |p T_xF|+|q T_xF|d\H^d(x)\\
&\le\int_F\ \lambda\ d\H^d(x)=\lambda \H^d(F)\end{array}\ee
since $|p T_xF|+|q T_xF|\le \lambda$. \qed
%
%
%
 
 As a corollary to Proposition 2.12, we have
 \begin{cor}Let $P_0^1,\cdots, P_0^m$ be $m$ mutually orthogonal planes of dimension $d$ in $\R^{md}$. Then their union $P_0=\cup_{i=1}^m P_0^i$ is minimal.
 \end{cor} 
 \nd By induction on $m$, with applying Proposition 2.12, and the fact that an $d$-plane is always a minimal set.\qed

Next we are going to deal with almost orthogonal cases:

\begin{lem}
 Let $P^i,1\le i\le m$, be $m$ planes of dimension $d$ in $\R^{md}$, and $\a^{ij}=(\a^{ij}_1,\cdots, \a^{ij}_d)$ be the characteristic angle between $P^i$ and $P^j$, $1\le i<j\le m$. Denote by $p^i$ the orthogonal projection on $P^i$. Then there exists $C_{m,d}(\a)$, with $\lim_{\a\to \frac\pi2}C_{m,d}(\a)=0$, such that for every unit simple $d-$vector $\zeta\in\wedge_d\R^{md}$, the sum of its projections to these $d$-planes satisfies
\be \sum_{i=1}^m |p^i(\zeta)|\le 1+C_{m,d}(\a),\ee
where $\a=\min_{1\le i<j\le m}\a_1^{ij}.$
\end{lem}

\nd We prove it by induction on $m$.

Case for $m=2$: Let $P^1$ and $P^2$ be two $d$-planes with characteristic angles $0\le \a_2\le\cdots\le \a_d\le\frac\pi2$. There there exists an orthonormal basis $\{e_1,\cdots e_{2d}\}$ of $\R^{2d}$ such that $P^1=P(\wedge_{i=1}^d e_i)$ and $P^2=P(\wedge_{i=1}^d(\cos\a_ie_i+\sin\a_ie_{i+d}))$.
Denote also by $p$ the orthogonal projection on ${P^1}^\perp=P(\wedge_{i=d+1}^{2d}e_i)$. Then 
\be |p^1(\zeta)|+|p^2(\zeta)|\le |p^1(\zeta)|+|p(\zeta)|+|(p^2-p)(\zeta)|\le 1+|(p^2-p)(\zeta)|.\ee

By Lemma \ref{F1} we know that
\be |p^1(\zeta)|+|p(\zeta)|\le 1.\ee
Estimate the last term, we have
\be \begin{split}|(p^2-p)(\zeta)|=&|<\wedge_{i=1}^d(\cos\a_ie_i+\sin\a_ie_{i+d})-\wedge_{i=d+1}^{2d}e_i,\zeta>|\\
=&|<\sum_{i=1}^{d+1}(\wedge_{j<i}\sin\a_je_j)\wedge\cos\a_ie_i\wedge(\wedge_{j>i}(\cos\a_je_j+\sin\a_je_{j+d})-\wedge_{i=d+1}^{2d}e_i,\zeta>|\\
\le&\sum_{i=1}^d|<(\wedge_{j<i}\sin\a_je_j)\wedge\cos\a_ie_i\wedge(\wedge_{j>i}(\cos\a_je_j+\sin\a_je_{j+d}),\zeta>|\\
&+|<\wedge_{i=1}^d\sin\a_ie_{j+d}-\wedge_{i=d+1}^{2d}e_i,\zeta>|\\
\le&\sum_{i=1}^d|\cos\a_i|+|1-\Pi_{i=1}^d\sin\a_i|\le d\cos\a_1+(1-\sin^2\a_1)\\
\le& d\cos\a_1+\cos^2\a_1\le(d+1)\cos\a_1.\end{split}\ee

Therefore
\be |p^1(\zeta)|+|p^2(\zeta)|\le 1+(d+1)\cos\a_1.\ee

Now suppose that (2.28) is true for $m-1$. Now denote by $P$ the $d-$plane $(\oplus_{1\le i\le m-1} P^i)^\perp$, $p$ the orthogonal projection on $P$, and $q$ the orthogonal projection on $(\oplus_{1\le i\le m-1} P^i)$. Then by Lemma \ref{F1},
\be |p(\zeta)|+|q(\zeta)|\le 1,\ee
and hence
\be\begin{split}
\sum_{i=1}^m|p^j(\zeta)|&=\sum_{i=1}^{m-1}|p^j\circ q(\zeta)|+|p^m(\zeta)|\\
&\le \sum_{i=1}^{m-1}|p^j\circ q(\zeta)|+|p(\zeta)|+|(p-p^m)(\zeta)|.\end{split}\ee

By the induction hypothesis,
\be \sum_{i=1}^{m-1}|p^j\circ q(\zeta)|\le (1+C_{m-1,d}(\a))|q(\zeta)|\ee
with $\lim_{\a\to\frac\pi2}C_{m-1,d}(\a)=0$,
therefore
\be\begin{split}
\sum_{i=1}^m|p^j(\zeta)|&\le (1+C_{m-1,d}(\a))|q(\zeta)|+|p(\zeta)|+|(p-p^m)(\zeta)|\\
&\le 1+C_{m-1,d}(\a)|q(\zeta)|+|(p-p^m)(\zeta)|\\
&\le 1+C_{m-1,d}(\a)+|(p-p^m)(\zeta)|.
   \end{split}
\ee

Now when $\a$ goes to $\frac\pi 2$, the angle between $P$ and $P^m$ goes to 0. Hence $|(p-p^m)(\zeta)|\le C'_{m,d}(\a)$ with $\lim_{\a\to\frac\pi2}C'_{m,d}(\a)=0$, and thus we get the conclusion.\qed

As a particular case of Lemma 2.27, if $P_0^1,\cdots P_0^m$ are $m$ mutually orthogonal $d$-planes in $\R^{md}$, then
\be \sum_{i=1}^m|p_0^i(\xi)|\le 1\ee
for all unit simple $d$-vector $\xi\in\wg_d\R^{dm}$.

Now set $\Xi(m,d):=\{\xi\in\wedge_d\R^{dm}\mbox{ unit simple },\sum_{i=1}^m|p_0^i(\xi)|=1\}$. For the purpose of next section, we want to decide $\Xi(m,d).$

\begin{lem}1) If $d\ge 3$, then
\be\Xi(m,d)=\{\xi\in\wedge_d\R^{dm}\mbox{ unit simple }, P(\xi)=P_0^i\mbox{ for some }1\le i \le m\};\ee
2) If $d=2$, 
\be\Xi(m,2)=\{(\sum_{i=1}^ma_iu_i)\wedge(\sum_{i=1}^ma_iv_i):v_i,u_i\in P_0^i\mbox{ unit, }v_i\perp u_i,a_i>0,\sum_ia_i^2=1\}.\ee
\end{lem}

\nd

1) $m\ge 3$: Denote by $P=P_0^m$ and $Q$ the space generated by $\cup_{i=1}^{m-1}$. Then $P$ and $Q$ are orthogonal and hence satisfiy all the hypothesis in Lemma \ref{F1}. In particular, since $d\ge 3$, (2.10) is true. Now for any unit simple vector $\xi$, by Lemma 2.27,
\be \sum_{i=1}^{m-1}|p_0^i(\xi)|=\sum_{i=1}^{m-1}|p_0^i\circ q(\xi)|\le |q(\xi)|.\ee

So if $\xi\in\Xi(m,d)$, then
\be 1=\sum_{i=1}^m|p_0^i(\xi)|=|p(\xi)|+\sum_{i=1}^{m-1}|p_0^i(\xi)|\le |p(\xi)|+|q(\xi)|\le 1.\ee
Hence $ |p(\xi)|+|q(\xi)|= 1$, and  by (2.11), 
\be \xi \mbox{ belongs to } P\mbox{ or }Q.\ee
By induction, we have $\xi$ belongs to one of the $P_0^i,1\le i\le m$, which yields (2.39).

2) $d=2$. Notice that in this case, the argument for $d\ge 3$ does not work, because (2.10) no longer holds. 

We prove (2.39) by induction en $m$. 

When $m=2$, this is just Wirtinger's inequality stated in 1.8.2 of \cite{Fe}, with $\nu=2$, $\R^4=\C_1\oplus\C_2$, $P_0^1=\C_1$, ${P_0^1}^\perp=\C_2$, $\mu=1$.

Suppose now that it is true for $m-1$.

Denote by $Q=\cup_{i=1}^{m-1}P_0^i$, and $q$ the projection on $Q$. If $x\wedge y\in \Xi$ with $x\perp y$, $x,y$ unit, then there exists $\theta_1,\theta_2\in[0,\frac\pi2]$ such that 
\be x=\cos\theta_1 p_0^m(x)+\sin\theta_1q(x),y=\cos\theta_2 p_0^m(y)+\sin\theta_2q(y).\ee

Hence $\sum_{i=1}^m|p_0^i(x\wedge y)|=1$ implies that $|q(x\wedge y)|+|p_0^m(x\wedge y)|=1$. By the same proof for the case $m=2$, we have
\be \theta_1=\theta_2,p_0^m(x)\perp p_0^m(y),q(x)\perp q(y), |q(x\wedge y)|=\sum_{i=1}^{m-1}|p_0^i[q(x)\wedge q(y)]|.\ee
By induction hypothesis we obtain the conclusion.\qed

\section{Uniqueness of $P_0$}

Now all we have to do is to prove the uniqueness of $P_0$. Recall that $P_0=\cup_\perp P_0^i$ is the union of $m$ orthogonal planes of dimension $d$ ($d\ge 2$). Denote by $p_0^i$ the orthogonal projection on $P_0^i,1\le i\le m$.

\begin{thm}[uniqueness of $P_0$]\label{unicite}Suppose that $E\subset \overline B(0,1)$ is a $d$-dimensional closed reduced set which is Almgren minimal in $B(0,1)\subset\ \R^{dm}$, and which satisfies that
\be p_0^i(E\cap B(0,1))\supset P_0^i\cap B(0,1), \forall  1\le i\le m;\ee
\be E\cap \partial B(0,1)=P_0\cap \partial B(0,1);\ee
\be \begin{split}H^d(E\cap &B(0,1))\le mv(d)\\
or\ equivalently\ &H^d(E\cap B(0,1))=mv(d).\end{split}\ee
where $v(d)=H^d(\R^d\cap B(0,1))$.

Then $E=P_0\cap \overline B(0,1)$.

\end{thm}

\nd

Take a set $E$ that satisfies all the hypotheses in the proposition.

First we denote still by $\Xi:=\{\xi\in\wedge_d\R^{dm}\mbox{ unit simple },\sum_{i=1}^m|p^i(\xi)|=1\}$. Then by the hypothesis of Proposition 3.1, all the inequalities in the hypotheses and in the proof of Lemma 2.15 with $\lambda=1$ are equalities. Hence we have

\begin{lem}
1) For $\H^d$-almost all $x\in E$, $T_xE\in P(\Xi)$.

2) For every $1\le i\le m$, for $\H^d$-almost all $z\in P_0^i\cap B(0,1)=p_0^i(E)$,  
\be N(p_0^i,z)=\sharp\{{p_0^i}^{-1}(z)\cap E\}=1.\ee 
\end{lem}

\bigskip

After the lemma, naturally we have to use (1), and hence to look at the set $P(\Xi)$. By Lemma 2.38, we have to prove the theorem for two cases: $d>2$ and $d=2$.

 \noindent\textbf{1st case: $d>2$.} We first prove it for $m=2$. 
 
By Lemma 2.38,
\be P(\Xi)=\{P_0^i,1\le i\le m\}.\ee
Then by the $C^1$ regularity (Proposition \ref{e1}), around every point $x\in E$ such that $T_xE$ exists, there exists $r_x>0$ and $1\le i\le n$ such that 
\be E\cap B(x,r_x)=(P_0^i+x)\cap B(x,r_x).\ee

Now we are going to deal with points that do not admit any tangent plane. Let $y\in E$ be such a point. Let $K$ be any blow-up limit of $E$ on $y$. We are going to prove that $K=P_0$ (and hence $P_0$ is the only blow-up limit of $E$ on $x$).

\begin{lem}\label{explosion2}For all $x\in K$ such that the tangent plane $T_xK$ of $K$ at $x$ exists, $T_xK=P_0^1$ or $P_0^2$.
\end{lem}

\nd Since $K$ is a cone, we can suppose that $|x|=1$.

K is a blow-up limit of $E$ on $y$, hence there exists a sequence $r_k$ such that $\lim_{ k\to \infty}r_k=0$ and
\be r_k^{-1}(E-y)\cap B(0,2)\stackrel{d_H}{\longrightarrow} K\cap B(0,2).\ee
Hence
\be r_k^{-1}(E-y)\cap B(x,r)\stackrel{d_H}{\longrightarrow} K\cap B(x,r)\ee
uniformly for $r\in(0,\frac12)$.

But the tangent plane $T_xK$ of $K$ on $x$ exists, hence $T_xK$ is the blow-up limit of $K$ on $x$, therefore there exists $0<r<\frac 12$ such that $d_H(K\cap B(x,r), P\cap B(x,r))< \frac12 r\e_1$, where $P=T_xK+x$, and $\e_1$ is as in Proposition \ref{e1}.

Fix this $r$, there exists $N>0$ such that for all $k>N,$ $d_H(r_k^{-1}(E-y)\cap B(x,r), K\cap B(x,r))<\frac 12 r\e_1$, because of (3.11). Hence 
\be d_H(r_k^{-1}(E-y)\cap B(x,r), P\cap B(x,r))<\e_1r.\ee

Denote by $E_k=r_k^{-1}(E-y)$, then $E_k$ is also minimal. Moreover
\be d_H(E_k\cap B(x,r), P\cap B(x,r))<\e_1r.\ee

Now after Proposition \ref{e1}, for $k$ large enough, $E_k$ is the graph of a $C^1$ map $f_k$ from $P$ to $P^\perp$ in $B(x,\frac34r)$.

Set $x_k=r_kx+y$, then $r_k^{-1}(x_k-y)=x$. Therefore
\be E_k\cap B(x,r)=r_k^{-1}((B(x_k,r_kr)\cap E)-y)\ee
 for $k>N$.
This means that in $B(x_k,\frac34 r_kr)$, $E$ is a $C^1$ graph on a plane $Q_k$ parallel to $T_xK$. However, for almost all $z\in B(x_k,\frac34 r_kr)\cap E$, the tangent plane of $E$ at $z$ exists and is $P_0^1$ or $P_0^2$, and the map $z\to T_zE$ is continuous on the $C^1$ graph, hence $E$ coincides with $P_0^1$ or $P_0^2$ in $B(x_k,\frac34r_kr)$. In other words, $E_k\cap B(x,\frac34r)=(P_0^1+x)\cap B(x,\frac34r)$ or $E_k\cap B(x,\frac34r)=(P_0^2+x)\cap B(x,\frac34r)$. Then by (3.11), in $B(x,\frac34r)$, $P$ is the limit of a sequence of planes, which are either $P_0^1+x$ either $P_0^2+x$. Therefore $P=P_0^1+x$ or $P_0^2+x$, hence $T_xK=P_0^1$ or $P_0^2$.\qed
 
\begin{lem}
 $K=P_0$. 
\end{lem}

\nd Denote by $K^i=\{x\in K:T_xK=P_0^i\}$, then we claim that 
\be K^i\subset P_0^i.\ee

In fact, if $x\in K^1, x\ne 0$, then since $K$ is a cone, we have $[0,x]\in K$. But $T_xK=P_0^1$, and for almost all $z\in K$ its tangent plane is $P_0^1$ or $P_0^2$, hence by an argument similar to that of Lemma 3.9, there exists a radius $r=r(x)>0$ such that $K\cap B(x,r)$ is a plane parallel to $P_0^1$. But $[0,x]\cap B(x,r)\subset K\cap B(x,r)$, hence $[0,x]\subset P_0^1$. In particular, $x\in P_0^1$. Hence we get
\be K^1\subset P_0^1,\ee
and similarly
\be K^2\subset P_0^2.\ee
But $K$ is minimal, therefore it is rectifiable, such that almost all point of $K$ admit a tangent plane. Then by Lemma 3.5, we have
\be H^m(K\bs (K^1\cup K^2))=0,\ee 
and hence
\be H^m(K\bs (P_0^1\cup P_0^2))=0,\ee
therefore
\be K\subset P_0,\ee
since $K$ is a reduced closed set.

Now if $K\ne P_0$, then there exists $x\in P_0\cap\partial B(0,1)$ such that $(0,x]\not\subset K$, because $K$ is a cone. Suppose for example that $x\in P_0^1$.

But $K$ is also closed, therefore there exists $r>0$ such that $B(\frac12x,r)\cap K=\emptyset$. In other words $K$ has a hole in the plane $P_0^1$. Thus we can easily deform $P_0^1\cap B(0,1)\bs B(\frac12 x,r)$ in $B(0,1)$ to a set of arbitrarily small measure, while fixing $\partial B(0,1)$ and $P_0^2$ at the same time. This implies that $\H^d(K\cap P_0^1)=0$, since $K$ is minimal. Hence $K=P_0^2$, which contradict the fact that $K$ is not a plane.

As a result, $K=P_0$.\qed

After the discussion above, in $E$ we have only two types of points : points of type $\P$, and points of type $P_0$. And for both types, the blow-up limit is unique.

Next we are going to give some regularity around a point $x$ of type $P_0$.

\begin{lem}
Let $x\in E$ be such that the blow-up of $E$ on $x$ is $P_0$. Then there exists $r>0$ such that $E\cap B(x,r)=(P_0+x)\cap B(x,r)$.
\end{lem}

\nd By the proofs of Lemmas 3.9 and 3.15, $P_0$ is the unique blow-up limit of $E$ on $x$. Hence there exists $r_0>0$ such that for all $r<r_0$, 
\be d_{x,r}(E,P_0+x)<\min\{\frac{1}{100},\frac{1}{10}\e_1\},\ee
where $\e_1$ is the one in Proposition \ref{e1}.

Denote by $C^i(x,s)={p_0^i}^{-1}(B(0,s)\cap P_0^i)+x$ (a ``cylinder"). 

For all $y\in E\cap B(x,\frac45r)\bs C^1(x,\frac15r)$, we have 
\be d_{y,\frac{1}{10}r}(E,P_0^1+x)<\e_1.\ee

In fact, by (3.23), (3.24) is true if we replace $P_0^1$ by $P_0$. On the other hand we know that $d(B(y,\frac{1}{10}r),P_0^2)>\frac{1}{10}r$ because $y\not\in C^1(x,\frac15r)$, hence we have (3.24).

Then by Proposition \ref{e1}, for each $y\in E\cap B(x,\frac45r)\bs C^1(x,\frac15r)$, in $B(y,\frac{3}{40}r)$, $E$ is a $C^1$ graph of $P_0^1$, and in particular $T_yE\ne P_0^2$. As a result $T_yE=P_0^1$. By (3.8),
 \be \mbox{there exists }r'_y>0\mbox{ such that }E\cap B(y,r'_y)=(P_0^1+y)\cap B(y,r'_y).\ee

Now fix a $y_r\in E\cap B(x,\frac45r)\bs C^1(x,\frac15r)$. Denote by $A=E\cap (P_0^1+y_r)\cap B(x,\frac45r)\bs C^1(x,\frac15r)$. Then $A$ is relatively closed in $(P_0^1+y_r)\cap B(x,\frac45r)\bs C^1(x,\frac15r)$, because $E$ is. And $A$ is non empty because $y_r\in A$. But $A$ is also open in $(P_0^1+y_r)\cap B(x,\frac45r)\bs C^1(x,\frac15r)$, because for each $y\in A$, (3.25) is true.

As a result $A=(P_0^1+y_r)\cap B(x,\frac45r)\bs C^1(x,\frac15r)$. Therefore we have
\be (P_0^1+y_r)\cap B(x,\frac45r)\bs C^1(x,\frac15r)\subset E\cap B(x,\frac45r)\bs C^1(x,\frac15r).\ee
But by (3.25), for each $y\in (P_0^1+y_r)\cap B(x,\frac45r)\bs C^1(x,\frac15r)$, 
\be E\mbox{  coincides with }P_0^1+y\mbox{ in }B(y,\frac{3}{40}r),\ee
hence
\be P_0^1+y=P_0^1+y_r,\ee
and
\be E\mbox{ coincides with }P_0^1+y_r\mbox{ in } B((P_0^1+y_r)\cap B(x,\frac45r)\bs C^1(x,\frac15r),\frac{3}{40}r).\ee

We know that
\be d_{B(x,r)\bs C^1(x,\frac{1}{10}r),r}(E,P_0^1+x)<\frac{1}{100},\ee
and
\be d(y_r,P_0^1+x)<\frac{1}{100}r,\ee
hence
\be d_{B(x,r)\bs C^1(x,\frac{1}{10}r),r}(E,P_0^1+y_r)<\frac{1}{20},\ee
which implies that 
\be \begin{split}E\cap B(x,\frac45r)\bs C^1(x,\frac15r)\subset B((P_0^1+y_r),\frac{3}{40}r)\cap B(x,\frac45r)\bs C^1(x,\frac15r)\\
\subset B((P_0^1+y_r)\cap B(x,\frac45r)\bs C^1(x,\frac15r),\frac{3}{40}r).\end{split}\ee
As a result, (3.29) implies that 
\be E\cap B(x,\frac45r)\bs C^1(x,\frac15r)=(P_0^1+y_r)\cap B(x,\frac45r)\bs C^1(x,\frac15r).\ee

Notice that (3.34) is true for all $r<r_0$. Hence for all $r<r_0$, $p_0^2(y_r)$ are the same, which is equal to  $p_0^2(x)$, because (3.23) is true for all $r$ small enough. Thus we have 
\be (P_0^1+x)\cap B(x,\frac45 r_0)\subset E\cap B(x,\frac45 r_0).\ee

Similarly we have
\be (P_0^2+x)\cap B(x,\frac45 r_0)\subset E\cap B(x,\frac45 r_0),\ee
and hence
\be (P_0+x)\cap B(x,\frac45 r_0)\subset E\cap B(x,\frac45 r_0).\ee

Then we claim that
\be (P_0+x)\cap B(x,\frac12 r_0)= E\cap B(x,\frac12 r_0).\ee
In fact, suppose $y\in E\cap  B(x,\frac12 r_0).$ Take $r=2|y-x|<r_0$, and suppose for example that $|p_0^1(y-x)|\ge|p_0^2(y-x)|$.  Then $y\in B(x,\frac45r)\bs C^1(x,\frac15r)$, which gives that $y\in P_0^1+x$ by (3.34).\qed

Now we get back to the proof of Theorem 3.1 for $d> 2$ and $m=2$.

Denote by $E^i=\{x\in E\cap B(0,1);T_xE\mbox{ exists and is }P_0^i\}$. Then for $i=1,2$, $E^i$ is open in $E$, by (3.8). In addition, at least one of $E^i$ is non empty, thanks to Lemma 3.5(1). Suppose for example that $E^1\ne\emptyset$. Denote also by $E^0$ the set of points of type $P_0$. Then $E$ is the disjoint union of $E^1,E^2$ and $E^0$. As a result $E^0\cup E^1=E\bs E^2$ is closed.

Let $x\in E^1$. We claim then
\be(P_0^1+x)\cap B(0,1)\subset E^1\cup E^0.\ee

In fact, denote by $A=(P_0^1+x)\cap B(0,1)\cap (E^1\cup E^0)$. Then $A$ is non empty, since $x\in A$. $A$ is also closed in $(P_0^1+x)\cap B(0,1)$ because $E^1\cup E^0$ is. On the other hand, by (3.8) and Lemma 3.22, $A$ is also open in $(P_0^1+x)\cap B(0,1)$. Now since $(P_0^1+x)\cap B(0,1)$ is connected, we have $A=(P_0^1+x)\cap B(0,1)$, therefore (3.39) is true.

As a result, 
\be (P_0^1+x)\cap\partial B(0,1)\subset E\cap \partial B(0,1).\ee
Then by (3.3), $p_0^2(x)=0$. 

Thus we obtain 
\be E^1\subset P_0^1\subset E^1\cup E^0\subset E,\ee 
so that $\H^d(E^1)\le \H^d(P_0^1\cap B(0,1))=v(d)$. But $\H^d(E^0)=0$. Hence by (3.4), $\H^d(E^2)>0$. In particular $E^2$ is non empty.

Similarly we have also 
\be E^2\subset P_0^2\subset E.\ee

Thus
\be P_0\subset E.\ee

We apply once again (3.4), and get $E=P_0$.

For the case $d>2$ and $m>2$, the proof is similar, so we will sketch it here:

Around a point that admit a tangent plane, (3.8) is always true; all other points in $E$ are of type $P_I:=\cup_{i\in I}P_0^i,I\subset \{1,2,\cdots,m\}$ (see the proof for Lemma 3.15). And if $x\in E$, $I\subset \{1,2,\cdots, m\}$ are such that $C_x=P_I$, where $C_x$ is the unique blow-up limit of $E$ on $x$, then there exists $r>0$ such that $E\cap B(x,r)=(P_I+x)\cap B(x,r)$. (see the proof for Lemma 3.22).

Denote by $E^i=\{x\in E,T_xE=P_0^i\}$, $E^0=E\bs\cup_iE^i$, and $E^{0i}=\{x\in E^0, P_0^i\subset C_x\}$, where $C_x$ denote the blow-up limit of $E$ on $x$, which verifies 2) above . Then $E^0=\cup_i E^{0i}$ is discrete, hence is finite. Every $E^i$ is open in $E$. Next we prove similarly that $E^i\cup E^{0i}=P_0^i$, and thus complete the proof for $d>2$ and $m>2$.

\bigskip

\noindent\textbf{Case 2: $d=2$.} For $d=2$ the set $\Xi$ is much larger (see Lemma 2.38). However we know more about regularity for 2-dimensional minimal sets. Note that if we identify each $P_0^i$ with a copy of the complex plane $\C$, and endow the space $\R^{2n}$ with the corresponding complex structure, then the set $\Xi$ is composed exactly of all holomorphic or anti-holomorphic planes.

By Lemma 3.5, $\H^2$-almost every point $x$ in $E$, $E$ admit a tangent plane $T_xE\in P(\Xi)$. Let us first look at those points that admit a tangent plane.

\begin{lem}
For each $1\le j\le m$, for every $x\in E$ such that $T_xE$ exists and is not $P_0^i$ for all $i\ne j$ (or equivalently, $|p_0^j(T_xE)|>0$), there exists $r_x>0$ such that in $B(x,r_x)$, $E$ coincides with the graph of a function $\varphi_x:P_0^j\to {P_0^j}^\perp=\oplus_{i\ne j}P_0^i$, where each function $\varphi_i=p_0^i\circ \varphi_x:P_0^j\to P_0^i$ is either analytic, or anti-analytic.
\end{lem}

\nd Since $T_xE\ne P_0^i$, by the $C^1$ regularity (Theorem \ref{c1}) of minimal sets and a direct application of implicit function theorem, the set $E$ is locally a $C^1$ graph of a function $\varphi_x:P_0^j\to {P_0^j}^\perp=\oplus_{i\ne j}P_0^i$. By Lemma 3.5 (1), the tangent plane of this graph at each point belong to $P(\Xi)$. Then the analyticity of anti-analyticity of the functions $\varphi_i=p_0^i\circ \varphi_x$ come from the structure of $P(\Xi)$, since all 2-vectors in $P(\Xi)$ are complex (See Lemma 2.38(2)). 

The readers could refer to Lemma 3.22 of \cite{2p} more detail, where the author gave the proof for $m=2$ in detail. \qed
%
%

\begin{lem}
There is no point of type $\Y$ in $E$.
\end{lem}

\nd 

Suppose that $x\in E$ is of type $\Y$. Denote by $C_x=\cup_{i=1}^3P_i$ its tangent cone, which is a $\Y$ set, where $P_i$ are three closed half planes that meet along a line $D$ which is generalized by a unit vector $v$. Denote by $Q_i$ the plane containing $P_i$. Then there exists unit vectors $w_i\in Q_i,1\le i\le 3$, $w_i\perp v$, and the  angle between any two of $w_i,1\le i\le 3$ is $120^\circ$. We want to show first that 
\be \mbox{at least one of those }Q_i\mbox{ does not belong to }P(\Xi).\ee

If $P_1\not\in P(\Xi)$, everything is fine. So suppose that $P_1\in P(\Xi)$. Then since $P_1=P(v\wedge w_1)$, by Lemmas 3.5(1) and 2.38 (2), there exists an orthonormal basis $\{e_j\}_{1\le j\le 2m}$ of $\R^{2m}$ with $P_0^i=P(e_i\wedge e_{m+i}$), and $a_i>0,\sum_{i=1}^ma_i^2=1,$ and $a_i=a_{i-m}$ for $i>m$, such that
\be v=\sum_{i=1}^m a_ie_i, w_1=\sum_{i=m+1}^{2m}a_ie_i.\ee

Then if we want $Q_2$ to be also in $P(\Xi)$, there should exist $\{\e_i\}_{n+1\le i\le 2m}$, each $\e_i$ equals to 1 or -1, such that
\be w_2=\sum_{i=m+1}^{2m}\e_ia_ie_i.\ee

Denote by $I=\{m+1\le i\le 2m,\e_i=1\}$, and $J=\{m+1\le i\le 2m,\e_i=-1\}$. Set $a_I=(\sum_{i\in I}a_i^2)^\frac12,a_J=(\sum_{i\in J}a_i^2)^\frac12$, $f_I=\frac{1}{a_I}\sum_{i\in I}a_ie_i$, $f_J=\frac{1}{a_J}\sum_{i\in J}a_ie_i$. Then $f_I,f_J$ are unit vectors mutually perpendicular, $a_I^2+a_J^2=1$, and 
\be w_1=a_If_I+a_Jf_j,w_2=a_If_I-a_Jf_J.\ee

If the angle between $w_1$ and $w_2$ is $120^\circ$, then they generate a 2-dimensional plane, and $a_I,a_J\ne 0$. 

Now if angles between $w_1,w_3$ and $w_2,w_3$ are also $120^\circ$, then these three $w_i$ belong to the same plane. As a result, $w_3\in P(f_I\wedge f_J)$, and the only possibility is that $a_I=\frac 12,a_J=\frac{\sqrt3}{2}$, and $w_3=-f_I$. But in this case, if we want $P(v\wedge w_3)\in P(\Xi)$, then for every $j\in J, a_j=0$. Thus $f_J=0$, and hence $w_1=w_2$. This is impossible because the angle between $w_1$ and $w_2$ is $120^\circ$.

Thus we get the claim (3.46). 

Now without loss of generality, suppose that $Q_1\not\in P(\Xi)$. Since $P(\Xi)$ is closed in $G(2m,2)$, and $Q_1\not\in P(\Xi)$, by the $C^1$ regularity of minimal sets around a $\Y$ points (Theorem \ref{c1} and Remark 1.30 (2)), there is a non empty relative open set $U\subset E$ such that for each $x\in U$, $T_xE$ exists but is not in $P(\Xi)$. But $U$ is open in $E$, hence it is of positive measure, this contradicts Lemma 3.5.\qed

\begin{lem}
For each $x\in E$ such that the tangent plane $E$ at $x$ does not exist, there exists $l$ planes $Q_1,Q_2,\cdots, Q_l\in P(\Xi)$, with $l\le m$ and $Q_i\perp Q_j$ for $i\ne j$, such that the blow-up limit $C_x$ of $E$ on $x$ is unique and is equal to $\cup_{i=1}^l Q_i$. 
\end{lem}

\nd 

Take any $x\in E$ and let $C_x$ be a blow-up limit of $E$ at $x$. Suppose also that $x=0$ for short. First we claim that
\be C_x\mbox{ doesn't contain any point of type }\Y.\ee

Suppose this is not true, then there exists $p\in C_x$ such that $p$ is of type $\Y$. Then $p$ is not the origin, because otherwise $C_x$ is of type $\Y$, and hence $0$ is of type $\Y$, which gives a contradiction with Lemma 3.45.

So $p$ is not the origin. Then since $C_x$ is a cone, for every $r>0$, $rp\in C_x$ is a point of type $\Y$. We can thus suppose that $||p||=1$. Then by our description of 2-dimensional minimal cones (Lemma 1.31), there exists $0<r<\frac12$ such that in $B(p,r)$, $C_x$ coincides with a cone $Y$ of type $\Y$ centered at $p$.

Define $d_{x,r}(E,F)=\frac 1r\max\{\sup\{d(y,F):y\in E\cap B(x,r)\},\sup\{d(y,E):y\in F\cap B(x,r)\}\}$, which is the relative distance of two sets $E,F$ with respect to the ball $B(x,r)$. Now $C_x$ is a blow-up limit of $E$ at $x$, so that there exists $s>0$ (large) such that $d_{0,2}(C_x,sE)<\frac{r\e_2}{100}$, where $\e_2$ is the constant in Proposition 16.24 of \cite{DJT} (this proposition says roughly that if a 2-dimensional minimal set is close enough to a $\Y$ set in $B(0,1)$, then it contains a $\Y$ type point in $B(0,\frac 12)$). Equivalently, $d_{p,\frac r2}(sE,C_x)<\frac{\e_2}{50}$.

It is not hard to show that \be d_{p,\frac r2}(sE,C_x)=d_{p,\frac r2}(sE,Y).\ee

Take a point $z\in sE$ such that $d(z,p)<\frac r2\times \frac{\e_2}{50}$, then $d_{z,\frac r4}(sE,Y+z-p)<\frac{\e_2}{10}$. Here $Y+z-p$ is a $\Y$ cone centered at $z$. But $sE$ is minimal (since $E$ is), therefore Proposition 16.24 of \cite{DJT} gives that $sE$ contains a $\Y$ point, and hence $E$, too. This contradicts Lemma 3.45. Thus we obtain our claim (3.51).

Since $C_x$ is a minimal cone, By Lemma 1.33, $C_x\cap\partial B(0,1)$ is a finite collection of great circles and arcs of great circles that meet by 3 with angles of $120^\circ$. Then (3.51) implies that there is no such arcs, since $C_x$ does not have $\Y$ points. As a result, $C_x\cap\partial B(0,1)$ is a finite collection of great circles, and therefore $C_x$ is the union of a finite number of transversal planes $\cup_{i=1}^l Q_i$, with $l\le m$. In particular, $C_x$ verifies the full-length property (Remark 1.32 (2)), hence by Theorem \ref{c1}, we have $C^1$ regularity around $x$, that is, there exists $r>0$ and a $C^1$ diffeomorphism $Phi:B(0,2r)\to \Phi(B(0,2r))$, such that $\Phi(0)=x, D\Phi(0)=Id$, and $E\cap B(x,r)=\Phi(X)\cap B(x,r)$. 

We claim hence that each $Q_i\in P(\Sigma)$. In fact, if this is not true, without loss of generality we suppose that $Q_1\not\in P(\Sigma)$. We know that the image $\Phi(Q_1)\cap B(x,r)\subset E$. Since $D\Phi(0)=Id$, the tangent plane to $\Phi(Q_1)$ is $Q_1$. But $P(\Sigma)$ is closed in $G(2m,2)$, and $Q_1\not\in P(\Sigma)$, $\Phi(Q_1)\cap B(x,r)$ is a $C^1$ manifold, hence there exists a neighborhood $U$ of $x$ such that for every $y\in U\cap \Phi(Q_1)$, the tangent plane of $\Phi(Q_1)$ on $y$ is not in $P(\Sigma)$. Note that  $U\cap \Phi(Q_1)\bs\{x\}$ is a subset of $E$ with positive measure, and $E$ admit a tangent plane at every $y\in U\cap \Phi(Q_1)\bs\{x\}$, which coincides with the tangent plane of $U\cap \Phi(Q_1)$ on $y$, hence $E$ contains a set of positive measure, such that on each point of this set, the tangent plane of $E$ is not in $P(\Sigma)$. This contradicts Lemma 3.5 (1).

Hence each $Q_i$ is in $P(\Sigma)$.

Next we prove $Q_i\perp Q_j$ for $i\ne j$. For example we prove it for $Q_1,Q_2$. 

Since $Q_1,Q_2\in P(\Sigma)$, by Lemmas 2.38 (2) , there exists unit vectors $v_k^1,u_k^1,v_k^2,u_k^2\in P_0^k$, $v_k^1\perp u_k^1,v_k^2\perp u_k^2,$ and $a_k^1,a_k^2>0$, $1\le k\le m$, such that $\sum_{k=1}^m{a_k^1}^2=\sum_{k=1}^m{a_k^2}^2=1$, and
\be Q_1=(\sum_{k=1}^ma_k^1v_k^1)\wedge (\sum_{k=1}^ma_k^1u_k^1),Q_2=(\sum_{k=1}^ma_k^2v_k^2)\wedge (\sum_{k=1}^ma_k^2u_k^2).\ee

So if $k$ is such that $a_k^1\ne 0$, then $p_0^k(Q_1)={a_k^1}^2>0$. After the  $C^1$ regularity of $E$ around $x$, we denote by $\varphi$ the local $C^1$ correspondence between $E$ and $C_x$ in a ball  $B(x,r)$. Then there exists $s>0$ such that $B(p_0^k(x),s)\cap P_0^k\subset p_0^k(\varphi(Q_1+x))$, since $Q_1$ is contained in the tangent cone of $E$ at $x$. This implies that $a_k^2=0$, because otherwise there exists $s'>0$ such that $B(p_0^k(x),s')\cap P_0^k\subset p_0^k(\varphi(Q_2+x)),$ and hence 
\be(B(p_0^k(x),\min\{x,s'\})\bs\{x\})\cap P_0^k\subset\{z\in P_0^k,\sharp\{{p_0^k}^{-1}\{z\}\cap E\}\ge 2\},\ee
which contradicts Lemma 3.5(2).

The above argument shows that for any $k$, $a_k^1$ and $a_k^2$ cannot be both non zero. This shows that $Q_1$ and $Q_2$ are orthogonal to each other. Thus the proof of Lemma 3.50 is completed. 
\qed

By Remark 1.32 (2), a union of transversal planes verifies the full length property if it is minimal. Hence by Lemma 3.50, and after Theorem \ref{c1}, for any point $x\in E$ which does not admit a tangent plane, we still have the $C^1$ regularity around it, in particular, $E$ admit a unique blow-up limit $C_x$ at $x$, which is tangent to $E$ at $x$, and which is a union of orthogonal planes in $P(\Xi)$. Thus we have

\begin{lem}
 For each $x\in E$, denote by $C_x$ a blow-up limit of $E$ on $x$, then there exists $r_x>0$ and a $C^1$ map $\varphi_x: B(x,2r_x)\to \varphi_x (B(x,2r_x))$ with $\varphi_x(x)=x,d\varphi(x)=Id$, and $E$ coincides with the graph of $\varphi_x$ of $C_x+x$ in $B(x,r_x)$.
\end{lem}

\begin{cor}
Let $x\in E$, denote by $C_x$ the blow up limit of $E$ on $x$ (and we know that this is the only blow-up limit because of the $C^1$ regularity). Suppose that $Q\subset C_x$ is a plane, and $i\in\{1,\cdots,m\}$ is such that $H^2(p_0^i(Q))>0$. Let $\varphi_x$ be as in Lemma 3.55. Then $\varphi_x(Q+x)$ is the graph of a $C^1$ function $\psi_x$ from $P_0^i$ to $\oplus_{j\ne i}P_0^i$, where for all $j$, $p_0^j\circ \psi_x$ is analytic or anti-analytic from $P_0^i$ to $P_0^j$.
\end{cor}

\nd The proof is the same as that of Lemma 3.44. \qed

By Lemmas 3.44, 3.50 and Corollary 3.56, we know that there exists countably many 2-dimensional $C^1$ manifolds $S_1,S_2\cdots S_n\cdots$, which are locally analytic or anti analytic, such that  
 $E\cap B(0,1)=\cup_i S_i$, and $S_i$ meet each other orthogonally. Then by the $C^1$ regularity, 
\be\overline S_l\bs S_l\subset \partial B(0,1).\ee

In fact if $x\in \overline S_l\cap B(0,1)\subset E\cap B(0,1)$, then there exists a $C^1$ ball $B(x,r_x)$ of $E$ on $x$, with $\varphi_x$ as the $C^1$ correspondence between $C_x$ and $E$. Hence there exists a plane $Q\subset C_x$ such that $Q=T_xS_l$. Therefore $\varphi(B(x,r_x)\cap (x+Q)\subset S_l$, and hence $x$ is an interior point of $\overline S_l$, which implies that $x\in S_l$.

\begin{lem}If there exists $x\in S_1$ such that the tangent plane $T_xS_1$ of $S_1$ on $x$ verifies $T_xS_1\perp P_0^1$, then $p_0^1(S_1)$ is a point. In other words, there exists $y\in P_0^1\cap B(0,1)$ such that $S_1\subset (y+{P_0^1}^\perp)$.
\end{lem}

\nd 

Let $x\in S_1$ be such that $T_xS_1\perp P_0^1$. Then there exists $i$ such that $|p_0^i(T_xS_1)|> 0$, and hence by Lemma 3.44 and Corollary 3.56, there exists $r_x>0$ such that in $B(x,r_x)$, $E$ coincides with the graph of a $C^1$ function $\varphi_x: P_0^i\to \oplus_{j\ne i}P_0^j$, and the function $\varphi_1=p_0^1\circ \varphi_x$ is analytic or anti-analytic, with $D\varphi_1=0$, i.e., of degree $d\ge 2$. So if $\varphi_1$ is not constant, then there exists $0<r<r_x$ such that $U=\varphi_1(B(x,r)\cap P_0^i)$ is an open set in $P_0^1$, and every point in $U\bs\{x\}$ has precisely $d$ pre-images. Hence the set $\{z\in P_0^1,\sharp \{{p_0^1}^{-1}(z)\cap E\}\ge 2\}$ is of positive measure, because it contains an open set $U$. This contradicts Lemma 3.5(2). 

Hence $\varphi_1$ is constant in $B(x,r_x)$. In other words, $p_0^1(S_1\cap B(x,r_x))=\{p_0^1(x)\}.$

Now denote by
\be A=\{y\in S_1,\mbox{ there exists }r_y>0\mbox{ such that }p_0^1(S_1\cap B(y,r_y))=\{p_0^1(y)\}\}.\ee

Then $A$ is open in $S_1$ and non empty. Moreover for each $y\in A$, $T_yS_1\perp P_0^1$. We want to prove that $A$ is also closed in $S_1$. For this purpose, take any sequence $\{y_l\}_{l\in\N}\subset A$ that converges to a point $y_0\in S_1$. Then by the $C^1$ regularity we have that $T_{y_0}S_1\perp P_0^1$, too. Then by the above argument, there exists $r>0$ such that $p_0^1(S_1\cap B(y_0,r))=\{p_0^1(y)\}.$ But $p_0^1(y)=\lim_{l\to\infty}p_0^1(y_l)=p_0^1(x)$, hence $y_0\in A$. Consequently $A$ is closed in $S_1$. Then since $S_1$ is connected, $A=S_1$, thus complete the proof of Lemma 3.58.\qed

As a direct corollary, we have

\begin{cor}If there exists $x\in S_1$ such that $T_xS_1$ is not perpendicular to $P_0^1$, then for all $y\in S_1$, $T_yS_1$ is not perpendicular to $P_0^1$.
\end{cor}

By (3.2), there exists $x\in E$ such that $|p_0^1(C_x)|>0$. Consequently there exists $l\in \N$ such that $x\in S_l$ and $T_xS_l$ is not perpendicular to $P_0^1$. Then by Corollary 3.60, 
\be \mbox{ for all }y\in S_l, T_yS_l \mbox{ is not perpendicular to }P_0^1.\ee

We claim that 
\be p_0^1(S_l)\supset P_0^1\cap B(0,1)\bs\{0\}.\ee

In fact, for each $z\in [P_0^1\cap B(0,1)\bs\{0\}] \cap p_0^1(S_l)$, by definition, there exists $y\in S_l$ such that $p_0^1(x)=z$. Then Lemma 3.44 and Corollary 3.56 implies that there exists $r_x>0$ such that $B(z,r_x)\cap P_0^1\subset p_0^1(S_l)$, and hence $p_0^1(S_l)$ is open in $P_0^1\cap B(0,1)\bs\{0\}$. On the other hand we are going to prove that  $p_0^1(S_l)$ is also closed in $P_0^1\cap B(0,1)\bs\{0\}$. 

So let $\{z_n\}_{n\in \N}\subset p_0^1(S_l)\cap P_0^1\cap B(0,1)\bs\{0\}$, which converges to a point $z_0\in P_0^1\cap B(0,1)\bs\{0\}$. Let $x_n\in S_l$ be such that $p_0^1(x_n)=z_n$. Then since $\overline S_l$ is compact, there exists a subsequence $\{x_{n_k}\}_{k\in\N}$ that converges to a point $x_0\in \overline S_l$. Then $p_0^1(x_0)=z_0$. If $x_0\not\in S_l$, by (3.57) we have $x_0\in\partial B(0,1)\cap  E=P_0\cap\partial B(0,1)$, which implies that $z_0\in\partial B(0,1)$ or $z_0=0$. But by hypothesis, $z_0\not\in\partial B(0.1)\cup\{0\}$, impossible. Hence $x_0\in S_l$, and hence $z_0\in \partial B(0.1)\cup\{0\}$. This proves that $p_0^1(S_l)$ is closed in $P_0^1\cap B(0,1)\bs\{0\}$, which contradicts the fact that $z_0\in B(0,1)$.

Hence $p_0^1(S_l)$ is both open and closed in $P_0^1\cap B(0,1)\bs\{0\}$, and obviously non empty. Then since $P_0^1\cap B(0,1)\bs\{0\}$ is connected, we get the claim  (3.62).

But $S_l\subset B(0,1)$, hence $p_0^1(S_l)\subset B(0,1)\cap P_0^1$. Thus we get
\be p_0^1(S_l)\bs\{0\}= P_0^1\cap B(0,1)\bs\{0\}.\ee

Next let us show that the projection $p_0^1$ is injective on $S_l$. Let $x_1,x_2\in S_l$ be such that $p_0^1(x_1)=p_0^1(x_2)=z\in P_0^1\cap B(0,1)$. Then there exists $r>0$ such that $B(x_1,r)\cap B(x_2.r)=\emptyset$. By Lemma 3.44 and Corollary 3.56, there exists $r'<r$ such that in $B(x_i,r')$, $E$ coincides with the graph of a $C^1$ function $\varphi_i$: $P_0^1\to {P_0^1}^\perp$. Consequently $B(z, r')\subset p_0^1(S_l\cap B(x_i,r')),i=1,2$. But $r'<r$ implies that $[S_l\cap B(x_1,r')]\cap [S_l\cap B(x_2,r')]=\emptyset$, therefore $\{z\in P_0^1:\sharp\{{p_0^1}^{-1}(z)\cap E\}\ge 2\}\supset B(z,r')\cap P_0^1$ is of positive measure. This contradicts Lemma 3.5 (2).

Thus we obtain the injectivity. As a result, $S_l\bs {p_0^1}^{-1}\{0\}$ is the graph of a $C^1$ function $\psi$ on $P_0^1\cap B(0,1)\bs \{0\}$, and moreover for each $i\ne 1$, the function $\psi_i=p_0^i\circ \psi$ is locally analytic or anti-analytic, and $|\psi_i(x)|<1$ for all $x\in P_0^1\cap B(0,1)\bs \{0\}$ (since the image is contained in $B(0,1)$. Then on applying Rad\'o's Theorem (cf. \cite{Ru} Thm 12.14), we know that $\psi$ is analytic or anti-analytic. (The readers could refer to the argument between (3.27) and Lemma 3.29 of \cite{2p} for more detail). Thus each $\psi_i$ is a bounded analytic or anti-analytic function on $P_0^1\cap B(0,1)\bs\{0\}$, therefore they can be extended to an analytic or anti-analytic function on $P_0^1\cap B(0,1)$ entirely, and $|\psi_i(0)|<1$, by the maximal principle. This implies that $(0,\psi(0))\in B(0,1)$. But $S_l$ is closed in $B(0,1)$, hence $(0,\psi(0))\in S_l$, this proves that in fact $S_l$ is a graph on $P_0^1\cap B(0,1)$, and $p_0^1(S_l)= P_0^1\cap B(0,1).$

Then the fact $S_l$ is a graph on $P_0^1\cap B(0,1)$ implies that $\overline S_l\bs S_l\subset {p_0^1}^{-1}[P_0^1\cap \partial B(0,1)]$. But $\overline S_l\bs S_l\subset E\cap\partial B(0,1)$, too, which gives
\be \overline S_l\bs S_l\subset P_0^1\cap \partial B(0,1).\ee
Therefore $\psi(z)$  tends to 0 when $z$ tends to $\partial B(0,1)$. Hence this is also true for $\psi_i$. Then by the maximum principle of analytic functions, every $\psi_i$ is constant, and equal to $0$.

Thus we deduce that $S_l=P_0^1\cap B(0,1)$. Consequently
\be P_0^1\cap B(0,1)\subset E.\ee

We treat similarly all the other $2\le i\le m$, and obtain that
\be P_0\cap B(0,1)=\cup_{1\le i\le m}P_0^i\cap B(0,1)\subset E.\ee

But by (3.4), $H^2(P_0\cap B(0,1))=m\pi=H^2(E)$, therefore
\be E\cap B(0,1)=P_0\cap B(0,1),\ee
where the conclusion of the uniqueness theorem follows.\qed

\section{A converging sequence of topological minimal competitors}

Now we begin to prove Theorem \ref{main}.

Fix any $d\ge 2,m\ge 2$.

For any $m$ $d-$planes $P^1,\cdots P^m$, with characteristic angles (see Definition 2.5) $\a^{ij}=(\a_1^{ij},\a_2^{ij},\cdots, \a_d^{ij})$ between $P^i$ and $P^j$, denote by $\a=(\a^{ij})_{1\le i<j\le m}$  the characteristic angle of these planes. For any $\theta\in[0,\frac\pi 2]$, then we say that the characteristic  angle of these $m$ $d-$planes $\a\ge\theta$ if for any $1\le i<j\le m$, $\min\a^{ij}\ge\theta$.

Now suppose that the conclusion of Theorem \ref{main} is not true. Then there exists a sequence of unions of $m$ $d$-planes of $P_k=\cup^{\a(k)}_{1\le i\le n} P^i_k\subset\R^{dm}$ with characteristic angle $\a(k)\ge\frac\pi2-\frac 1k$, which are not topologically minimal. Recall also that $P_0=\cup_{1\le i\le m} P_0^i$ is the orthogonal union of two planes. Choose an orthonormal basis $\{e_i\}_{1\le i\le md}$ of $\R^{dm}$ such that $P_0^i=e_{d(i-1)+1}\wg e_{d(i-1)+2}\wg\cdots \wg e_{di}, 1\le i\le m$. After necessary rotations, we suppose also that all the $P_k^1,k\ge 0$ are the same, and $d_H(P_0\cap \overline B(0,1), P_k\cap\overline B(0,1)\to 0$ when $k\to\infty$.

For each $k$, since $P_k$ is not topologically minimal, by definition of topological minimal sets, and the fact that $P_k$ is a cone, there exists a $F\in\mathscr F_k$ such that
\be \inf_{F\in\mathscr F_k}\H^d(F\cap B(0,1))<\H^d(P_k\cap B(0,1))=mv(d),\ee
where $v(d)$ denotes the $d$-Hausdorff measure of the unit ball in $\R^d$, and $\mathscr F_k$ is the class of all the topological competitors in $B(0,1)$. The following proposition will guarantee the existence of a topological competitor $F_k$ of $P_k$ in $B(0,1)$, which is topologically minimal in $B(0,1)$. 

%
%

\begin{pro}Let $P^1,\cdots P^m$ be $m$ mutually transversal planes of dimension $d$ in $\R^{dm}$. Denote by $B=B(0,1)$ the unit ball. Denote by $\mathfrak F$ the set of all the topological competitors of $E=\cup_{1\le i\le m}P^i$ in B (cf. Definition 1.16). Then there exists $F_0\in\mathfrak F$ such that 
\be \H^d(F_0\cap B)=\inf_{F\in\mathfrak F}\H^d(F\cap B).\ee
Moreover, $F_0\cap\overline B$ is contained in the convex hull of $E\cap \overline B$, and for any $1\le i\le m$, if we denote by $p^i$ the orthogonal projection from $\R^{dm}$ to $P^i$, then
\be p^i(F_0\cap B)\supset P^i\cap B.\ee
\end{pro}

\nd This is an easy corollary of the following theorem.

\begin{thm}[cf.\cite{topo} Theorem 4.28]Let $E\subset \R^n$ be a closed set. $B\subset \R^n$ is an open ball. Let $\{w_j\}_{j\in J}$ be a family of smooth $n-d-1$-surfaces in $\R^n\bs(\overline B\cup E)$, which are non-zero in $H_{n-d-1}(\R^n\bs E)$. Set
\be \mathfrak F=\{F\subset \R^n,F\bs B=E\bs B\mbox{ and for all }j\in J, w_j\mbox{ is non-zero in }H_{n-d-1}(\R^n\bs F)\}.\ee

Then there exists $F_0\in\mathfrak F$ such that
\be \H^d(F_0\cap B)=\inf\{\H^d(F\cap B);F\in\mathfrak F\}.\ee
Moreover, $F_0\cap B$ is contained in the convex hull of $E\cap B$.
\end{thm}

We apply Theorem 4.5, on taking $E$ to be $\cup _{1\le i\le m}P^i$ and $\{w_j\}_{j\in J}$ to be the family of all circles outside $B$ which are non-zero in $H_{md-d-1}(\R^{md}\bs \cup _{1\le i\le m}P^i)$, we get the existence of a set $F_0$ in $\mathfrak F$, and $F_0$ being contained in the convex hull of $\cup _{1\le i\le m}P^i\cap B$. The projection property (4.4) comes directly from the fact that $F_0$ is a topological competitor of $\cup _{1\le i\le m}P^i$.

\begin{rem}In \cite{2p}, things are probably more complicated, because we do not have any existence theorem for Almgren minimal sets there. Fortunately for topological minimal sets, we have the above existence theorem.
\end{rem}

Now by applying Proposition 4.2, we have the existence of a topologially minimal competitor $F_k$ of of $P_k$. Denote by $E_k=F_k\cap\overline B$ the part of $F_k$ inside $\overline B$. Then $E_k$ is contained in the convex hull $C_k$ of $P_k\cap\overline B$.  
\be p_k^i(E_k)\supset P_k^i\cap B,\ee
and  since $P_k$ is not minimal, by (4.1) and (4.7),
\be \H^d(E_k)<\H^d(P_k\cap B)=mv(d).\ee

Now since $\overline B$ is compact, we can extract a converging subsequence of $\{E_k\}$, denoted still by $\{E_k\}$ for short. Denote by $E_\infty$ their limit. Then $E_\infty$ is contained in $\cap_n\cup_{k>n} C_k$, such that $E_\infty\cap \partial B\subset (\cap_n\cup_{k>n} C_k)\cap\partial B=P_0\cap\partial B$. On the other hand $E_\infty\cap \partial B\supset\lim_{k\to\infty} (P_k\cap\partial B)=P_0\cap\partial B$. Hence 
\be E_\infty\cap\partial B=P_0\cap\partial B.\ee

We want to use the uniqueness theorem \ref{unicite}, to prove that $E_\infty$ is in fact $P_0\cap B$. So we have to check all the conditions:

1) First, we know that $E_k$ are topologically minimal in $B$, and hence are Almgren minimal in $B$. Hence as their limit, $E_\infty$ is Almgren minimal in $B$;

2) Since $E_\infty$ is the limit of $E_k$, and $P_0$ is the limit of $P_k$, hence the projection property (3.2) comes from (4.9);

3) The boundary condition (3.3) is already proved in (4.11); the measure estimate (3.4) is guaranteed by (4.10).

Recall that $E_k=F_k\cap \overline B$, and $F_k\bs B=P_k\bs B$, where $P_k$ converge to $P_0$. Thus we have a sequence of closed sets $F_k$, each $F_k$ is a minimal topological competitor of $P_k$ in $B$, and $F_k$ converge to $P_0$.

\section{A stopping time argument}

We will continue our argument, by cutting each $E_k$ into two pieces. One piece is inside a small ball near the origin, where something complicated happens there, and we can only estimate its measure by projection argument; the other piece is outside the small ball, where $E_k$ is very near $P_k$, and by the regularity of minimal sets near planes, $E_k$ is composed of $m$ $C^1$ graphs on $P_k^i, 1\le i\le m$ respectively, where we will estimate their measures, by harmonic extensions. So the first step is to find this small ball, with the critical radius, by a stopping time argument.

For each $k$ and $i=1,\cdots, m$, denote by 
\be C_k^i(x,r)={p_k^i}^{-1}(B(0,r)\cap P_k^i)+x,\ee
and
\be D_k(x,r)=\cap _{1\le i\le n}C_k^i(x,r).\ee
Notice that $D_k(0,1)\supset B(0,1)$ and $D_k(0,1)\cap P_k=B(0,1)\cap P_k$.

We say that two sets $E,F$ are $\e r$ near each other in an open set $U$ if
\be d_{r,U}(E,F)<\e,\ee 
where
\be d_{r,U}(E,F)=\frac 1r\max\{\sup\{d(y,F):y\in E\cap U\},\sup\{d(y,E):y\in F\cap U\}\}.\ee
We define also
\be\begin{split}
    &d_{x,r}^k(E,F)=d_{r,D_k(x,r)}(E,F)\\
&=\frac 1r\max\{\sup\{d(y,F):y\in E\cap D_k(x,r)\},\sup\{d(y,E):y\in F\cap D_k(x,r)\}\}.
   \end{split}\ee

\begin{rem}
 Observe that $d_{r,U}(E,F)\ne\frac 1rd_H(E\cap U,F\cap U)$. For example, take $U=B(0,1)$, set $E_n=\partial B(0,1-\frac 1n)$, and $F_n=\partial B(0,1+\frac 1n)$, then we have \be d_{1,U}(E_n,F_n)\to 0, d_H(E_n\cap U,F_n\cap U)=d_H(E_n\cap U,\emptyset)=\infty.\ee
\end{rem}

Now we start our stopping time argument. We fix a $\e$ small and a $k$ large, and we set $s_i=2^{-i}$ for $i\ge 0$. Denote by $D(x,r)=D_k(x,r), d_{x,r}=d^k_{x,r}$ for short. Then we proceed as follows.

Step 1: Denote by $q_0=q_1=0$, then in $D(q_0,s_0)$, the set $E_k$ is $\e s_0$ near $P_k+q_1$ when $k$ is large, because $E_k\to P_0$ and $P_k\to P_0$ implies that $d_{0,1}(E_k,P_k)\to 0$.

Step 2: If in $D(q_1,s_1)$, there is no point $q\in\R^{dm}$ such that $E_k$ is $\e s_1$ near $P_k+q$, we stop here; otherwise, there exists a point $q_2$ such that $E_k$ is $\e s_1$ near $P_k+q_2$ in $D(q_1,s_1)$. Here we ask $\e$ to be small enough (say, $\e<\frac {1}{100}$) such that such a $q_2$ is automatically in $D(q_,\frac12 s_1)$, by the conclusion of the step 1. Then in $D(q_1,s_1)$ we have simultaneously  
\be d_{q_1,s_1}(E_k,P_k+q_1)\le s_1^{-1}d_{q_0,s_0}(E_k,P_k+q_1)\le 2\e;\ d_{q_1,s_1}(E_k, P_k+q_2)\le\e.\ee
This implies that $d_{q_1,\frac12 s_1}(P_k+q_1,P_k+q_2)\le 12\e$ when $\e$ is small. And hence $d(q_1,q_2)\le 6\e$.

Now we are going to define our iteration process. Notice that this process depends on $\e$, hence we also call it a $\e-$process.

Suppose that $\{q_i\}$ are defined for all $i\le n$, with 
\be d_{q_1,q_{i+1}}\le 12 s_i\e=12\times 2^{-i}\e\ee
for $0\le i\le n-1$, and hence
\be d_{q_i,q_j}\le 24\e_{\min(i,j)}=2^{-\min(i,j)}\times 24\e\ee
for $0\le i,j\le n$, and that for all $i\le n-1$, $E_k$ is $\e s_i$ near $P_k+q_{i+1}$ in $D(q_i,s_i)$. We say in this case that the process does not stop at step $n$. Then

Step $n+1$: We look inside $D(q_n,s_n)$.

If $E_k$ is not $\e s_n$ near any $P_k+q$ in this ``ball'' of radius $s_n$, we stop. In this case, since $d(q_{n-1},q_n)\le 12\e s_{n-1}$, we have $D(q_n,2s_n(1-12\e))=D(q_n,s_{n-1}(1-12\e))\subset D(q_{n-1},s_{n-1})$, and hence
\be\begin{split}
    d_{q_n,2s_n(1-12\e)}(P_k+q_n,E_k)&\le (1-12\e)^{-1}d_{q_{n-1},s_{n-1}}(P_k+q,E_k)\\
&\le\frac{\e}{1-12\e}.
   \end{split}
\ee
Moreover
\be d(q_n,0)=d(q_n,q_1)\le 2^{-\min(1,n)}\times 24\e=12\e.\ee

Otherwise, we can find a $q_{n+1}\in\R^{dm}$ such that $E_k$ is still $\e s_n$ near $P_k+q_{n+1}$ in $D(q_n,s_n)$, then since $\e$ is small, as before we have $d(q_{n+1},q_n)\le 12\e s_n$, and for $i\le n-1$,
\be d_{q_i,q_n}\le\sum_{j=i}^n 12\times 2^{-j}\e\le 2^{-\min(i,n)}\times 24\e.\ee
Thus we get our $q_{n+1}$, and say that the process does not stop at step $n+1$.

We will see in the next section, that the process has to stop at a finite step. And for each $k$, if the process stop at step $n$, we define $o_k=q_n,r_k=s_n$. Then $D_k(o_k,r_k)$ is the critical ball that we look for, because inside the small ball, by definition we know that $E_k$ is $\e s_n$ far from any translation of $P_k$, but outside it, things are near. We also have, by (5.12), $d(o_k,0)\le 12\e$, hence the center $o_k$ of the critical ball is near the origin.

\section{Regularity and projection properties of $E_k$}

\begin{pro}\label{reg}
 There exists $\e_0\in(0,\frac{1}{100})$, such that for any $\e<\e_0$ fixed and for $k$ large, if our $\e-$process does not stop before the step $n$, then

(1) The set $E_k\cap (D_k(0,\frac{39}{40})\bs D_k(q_n,\frac{1}{10}s_n))$ is composed of $m$ disjoint pieces $G^i,i=1,\cdots, m$, such that 
\be G^i\mbox{ is the graph of a }C^1\mbox{ map }g^i:D_k(0,\frac{39}{40})\bs D_k(q_n,\frac{1}{10}s_n)\cap P_k^i\to {P_k^i}^\perp\ee
with
\be ||\nabla g^i||_\infty<1;\ee

(2) For each $t\in[\frac{1}{10}s_n,s_n],$
\be E_k\cap (D_k(0,1)\bs D_k(q_n,t))=G_t^1\cup G_t^2\cup\cdots\cup G_t^m,\ee
where the $G_t^i,1\le i\le m$ do not meet. Moreover
\be P_k^i\cap (D_k(0,1)\bs C_k^i(q_n,t))\subset p_k^i(G_t^i)\mbox{ for }i=1,\cdots, m,\ee
where $p_k^i$ is the orthogonal projection on $P_k^i,i=1,\cdots, m$;

(3) The projections $p_k^i:E_k\cap \overline D_k(q_n,t)\to P_k^i\cap \overline C_k^i(q_n,t),i=1,\cdots, m$ are surjective, for all $t\in[\frac{1}{10}s_n,s_n].$
\end{pro}

Before we give the proof, first we give a direct corollary of (2), which shows that the $\e$-process stated in the previous subsection has to stop at a finite step for any $k$.

\begin{cor}
 For any $k$ and $\e<\e_0$, the $\e$ process has to stop at a finite step.
\end{cor}

\nd Since $\H^d(E_k\cap \overline B(0,1))<m v(d)$, there exists $n_k>0$ such that
\be \inf_{q\in\R^{dm}} \H^d(P_k\cap \overline B(0,1)\bs D(q,s_{n_k}))>\H^d(E_k).\ee
Then our process need to stop before the step $n_k$, because otherwise, we use the term (2) in Proposition \ref{reg}, for $t=s_n$, and get the disjoint decomposition
\be E_k=[E_k\cap D(q_{n_k},s_{n_k})]\cup G_{s_{n_k}}^1\cup G_{s_{n_k}}^2,\ee
therefore
\be \begin{split}
     \H^d(E_k)&\ge \H^d(G_{s_{n_k}}^1)+\H^d(G_{s_{n_k}}^2)\ge \H^d[p_k^1(G_{s_{n_k}}^1)]+\H^d[p_k^2(G_{s_{n_k}}^2)]\\
&\ge \H^d(P_k\cap \overline B(0,1)\bs D(q_{nk},s_{n_k})>\H^d(E_k),
    \end{split}
\ee
which leads to a contradiction.\qed

Now we are going to prove Proposition \ref{reg}.

\noindent Proof of Proposition \ref{reg}.

For (1), notice that every topological minimal set is an Almgren minimal sets, hence (1) and (6.4) are direct corollaries of the proposition 6.1 (1) of \cite{2p}. 

As a result of (1), we know that (6.5) in (2) is true if we replace all the $D_k(0,1)$ with $D_k(0,\frac{39}{40})$. Hence we have to prove that 
\be P_k^i\cap D(0,1)\bs D(0,\frac{39}{40})\subset p_k^i(G_t^i).\ee

We prove it for $i=1$ for example. The other case is the same. 

We know that $G_t^i$ is very close to $P_k^i$, $2\le i\le m$, and when $k$ is large, $P_k^1$ and $P_k^i$ are almost orthogonal, hence the projection of $P_k^i\cap D(0,1)$ under $p_k^1$ is far away from $P_k^1\cap D(0,1)\bs D(0,\frac{39}{40}).$ On the other hand, the projection of the part $E_k\cap D(q_n,t)$ under $p_k^1$ is always contained in  $D(q_n,t)$, hence also far away from $P_k^1\cap D(0,1)\bs D(0,\frac{39}{40}).$ So (6.10) is equivalent to say that
\be P_k^1\cap D(0,1)\bs D(0,\frac{39}{40})\subset p_k^1(E_k\cap\overline D(0,1)).\ee

We are going to prove a stronger one, that is
\be P_k^1\cap \overline B(0,1)\subset p_k^1(E_k\cap\overline B(0,1)).\ee

So suppose that (6.12) is not true. That is, there exists $x\in P_k^1\cap \overline B(0,1)$ such that ${p_k^1}^{-1}\{x\}\cap (E_k\cap\overline B(0,1))=\emptyset.$ In other words, ${p_k^1}^{-1}\{x\}\cap\overline B(0,1)=(P+x)\cap\overline B(0,1)$ does not meet $E_k$, where by $P$ we denote the $md-d$-subspace orthogonal to $P_k^1$. But $E_k$ is closed, so there exists $\d>0$ such that the neighborhood $B((P+x)\cap\overline B(0,1),\d)$ does not meet $E_k$. Denote by $S=(P+x)\cap\partial B(0,1+\d)$, then $S$ is a $md-d-1$-sphere that does not meet $E_k$, and $S$ is zero in $H_{md-d-1}(\R^dm\bs E_k)$, because it is the boundary of the disc $(P+x)\cap\partial B(0,1)\subset\R^{md}\bs E_k$.

However, $S$ is non zero in $H_{md-d-1}(\R^{md}\bs P_k^1)$, hence is non zero in $H_{md-d-1}(\R^{md}\bs P_k)$. This contradicts the fact that $E_k$ is a topological competitor of $P_k$. 

We have thus (6.12), which gives (6.11), and hence (6.10). 

So we get (2).

To prove (3), the idea is almost the same as above. We prove it for $i=1$ for example, denote still by $P$ the orthogonal $md-d$-subspace of $P_k^1$. Suppose that there exists $x\in P_k^1\cap \overline C_k^1(q_n,t)$ satisfying ${p_k^1}^{-1}(x)\cap\overline D_k(q_n,t)\cap E_k=\emptyset$. Equivalently, 
\be (P+x)\cap\overline D_k(q_n,t)\cap E_k=\emptyset.\ee
Denote by $S=(P+x)\cap\partial D_k(q_n,t)$. Then $S$ is zero in $H_{md-d-1}(\R^{md}\bs E_k)$, since it is the boundary of the disc $(P+x)\cap\overline D_k(q_n,t)$.
But since the $\e$-process does not stop at step $n$, outside $D_k(q_n,t)$, $E_k$ is composed of two disjoint pieces that are $\e$ closed to $P_k^1$ and $P_k^2$ respectively, and $d(q_n,0)\le 12\e$, so in fact we can deform our $S$ to any circle $S_y=(P+y)\cap \partial B(y,\frac 12)$ for all $y\in P_k^1\bs\overline B(0,1)$. Hence such a $S_y$ is zero in $H_{md-d-1}(\R^{md}\bs E_k)$. Again, we know that such a $S_y$ is non zero in $H_{md-d-1}(\R^{md}\bs E_k)$, which contradicts the fact that $E_k$ is a topological competitor of $P_k$.

Thus we get (3). And the proof of Proposition 6.1 is finished.\qed

\section{Estimates for graphs by harmonic extension}

In this section we will give some estimates on the Dirichlet energy of a function with prescribed partial boundary condition. This will be used in the next section to estimate the measure of the graphs $G^i,1\le i\le m$ in (6.2), for each $E_k$. In the proof of the following propositions we will use the space of spherical harmonics. We will state some necessary definitions and theorems here. Please refer to \cite{St71} for more detail.

Given an integer $d\ge 2$, set 
\be \mathscr H_n=\mathscr H_n(\R^d)=\{p|_{S^{d-1}}:p\mbox{ is a homogeneous harmonic polynomials of degree }n\mbox{ in }\R^d\},\ee
 the space of spherical harmonics of degree $n$. Then we have the following properties:

\begin{pro}
$1^\circ$  Each $f\in \mathscr H_n$ is an eigen function of $\Delta_{S^{d-1}}$ associated to the eigenvalue $\lambda_n=-n(n+d-2)$, where $\Delta_{S^{d-1}}$ is the Laplace-Beltrami operator on $S^{d-1}$, the angular part of the Laplacian; 

$2^\circ$ The collection of all finite linear combinations of elements of $\cup_{n=0}^\infty\mathscr H_n$ is dense in $L^2{S^{d-1}}$;

$3^\circ$ Let $Y^{(n)}$ and $Y^{(m)}$ be spherical harmonics of degree $n$ and $m$, and $n\ne m$, then 
\be\int_{S^{n-1}}Y^{(n)}(\theta)Y^{(m)}(\theta)d\theta=0.\ee
\end{pro}

\begin{rem} By $1^\circ$, if $f\in\mathscr H_n$, then $r^{2-d-n}f(\theta)$ is also a harmonic function.
\end{rem}

Now if we consider $\mathscr H_n$ as a subspace of the Hilbert space $L^2(S^{d-1})$ with the scalar product $(f,g)=\int_{S^{d-1}}f(\theta)\overline{g(\theta)}d\theta$, we have the following corollary:

\begin{cor}\label{base}Denote by $a_n=dim\mathscr H_n<\infty$, and let $\{Y_1^{(n)},\cdots,Y_{a_n}^{(n)}\}$ be an orthonormal basis of $\mathscr H_n$, then $\cup_{n=0}^\infty \{Y_1^{(n)},\cdots,Y_{a_n}^{(n)}\}$ is an orthonormal basis of $L^2(S^{d-1})$. Moreover, for each $f\in L^2(S^{d-1})$, there exists a unique representation:
\be f=\sum_{n=0}^\infty\sum_{i=1}^{a_n}b_i^{(n)}Y_i^{(n)},\ee
where the series converges to $f$ under the $L^2$ norm. Hence we have
\be||f||^2_{L^2(S^{d-1})}=\sum_{n=0}^\infty\sum_{i=1}^{a_n}|b_i^{(n)}|^2.\ee
\end{cor}

Now we can start to give our estimates.

\begin{pro}Let $d\ge 2$ be an integer, $0<r_0<\frac12$ and $u_0\in C^1(\partial B(0,r_0)\cap \R^d,\R)$. Denote by $m(u_0)=(r_0^{d-1}s_{d-1})^{-1}\int_{\partial B(0,r_0)}u_0$ its average, where $s_{d-1}=H^{d-1}(S^{d-1})$, the $d-1$-measure of the unit sphere of $\R^d$.
Then for all $u\in C^1((\overline{B(0,1)}\backslash B(0,r_0))\cap \R^d, \R)$ that satisfies
\be u|_{\partial B(0,r_0)}=u_0,\ee
we have
\be\int_{B(0,1)\backslash B(0,r_0)}|\nabla u|^2\ge \frac13r_0^{-1}\int_{\partial B(0,r_0)}|u_0-m(u_0)|^2.\ee
\end{pro}

\nd

Let $u$ be a $C^1$ function as in the statement of the proposition. Then if $v$ is a solution of the equation 
\be\left\{\begin{array}{l}
\Delta v=0;\\
v|_{\partial B(0,r_0)}=u_0;\\
\frac{\partial u}{\partial \vec n}=0\mbox{ on }\partial B(0,1)    
   \end{array}\right.\ee
where $\vec n$ is the unit exterior normal vector on $\partial B(0,1)$, then $v$ minimizes Dirichlet's energy among all $C^1$ function $u$ on $\overline B(0,1)\bs B(0,r_0)$ with boundary condition $u|_{\partial B(0,r_0)}=u_0$.

So we are just going to look for a solution $v$ of the equation (7.11), and then prove (7.10) for $v$.

Let $\{Y_1^{(n)},\cdots,Y_{a_n}^{(n)}\}$ be an orthonormal basis of $\mathscr H_n(\R^d)$. We express $u_0$ under this basis
\be u(r_0,\theta)=u_0(\theta)=m(u_0)+\sum_{n=1}^\infty\sum_{i=1}^{a_n}B_i^{(n)}Y_i^{(n)}.\ee
Set
\be \begin{split}v&:\overline B(0,1)\backslash B(0,r_0)\to\R,\\
 v(r,\theta)&=m(u_0)+\sum_{n=1}^\infty(A_nr^n+B_nr^{2-d-n})\sum_{i=1}^{a_n}B_i^{(n)}Y_i^{(n)}.
\end{split}\ee
By Remark 7.4, $v$ is harmonic. We want $v$ to verify (7.11). Notice that
\be\begin{split} v&|_{\partial B(0,r_0)}=u_0\Leftrightarrow \\
&(A_nr_0^n+B_nr_0^{2-d-n})=1,\mbox{ for each }n\mbox{ and }1\le i\le a_n,
   \end{split}\ee
and
\be\begin{split} \frac{\partial v}{\partial \vec n}&=\sum_{n=1}^\infty(A_nnr^{n-1}+B_n(2-d-n)r^{1-d-n})\sum_{i=1}^{a_n}B_i^{(n)}Y_i^{(n)}=0\\
&\mbox{ on }\partial B(0,1)=\{(r,\theta):r=1\}.
   \end{split}\ee
So we ask
\be \left(\begin{array}{cc}r_0^n & r_0^{2-d-n}\\n & 2-d-n\end{array}\right)\left(\begin{array}{c}A_n\\B_n\end{array}\right)=\left(\begin{array}{c}1\\0\end{array}\right).\ee

The determinant of the coefficient matrix is $(2-d-n)r_0^n-nr_0^{2-d-n}$, with $n>0,r_0^{2-d-n}>r_0^n>0,2-d-n\le 0$, hence it is always strictly negative, therefore (7.16) admits always a solution
\be A_n=\frac{n+d-2}{(n+d-2)r_0^n+nr_0^{2-d-n}},B_n=\frac{n}{(n+d-2)r_0^n+nr_0^{2-d-n}}.\ee

Now let us calculate $\int_{\overline B(0,1)\bs B(0,r_0)}|\nabla v|^2$. To estimate $\nabla v$ we have
\be\nabla_{S^{d-1}}v(r,\theta)=\sum_{n=1}^\infty(A_nr^n+B_nr^{2-d-n})\sum_{i=1}^{a_n}B_i^{(n)}\nabla_{S^{d-1}}Y_i^{(n)},\ee
\be\frac{\partial v}{\partial r}=\sum_{n=1}^\infty (nA_nr^{n-1}+(2-d-n)B_nr^{1-d-n})\sum_{i=1}^{a_n}B_i^{(n)}Y_i^{(n)},\ee
and
\be
\int_{S^{d-1}}|\nabla v|^2d\theta=\int_{S^{d-1}}|\frac{\partial v}{\partial r}|^2+|\frac 1r\nabla_{S^{d-1}}v|^2d\theta.\ee
Notice that
\be <Y_i^{(n)},Y_j^{(m)}>_{L^2(S^{d-1})}=\delta_{(i,n)(j,m)}\mbox{ for }1\le n\mbox{ and }1\le i\le a_n,\ee
 where $\delta$ is the Kronecker symbol, and hence
\be\begin{split}<\nabla_{S^{d-1}} Y_i^{(n)},&\nabla_{S^{d-1}} Y_j^{(m)}>_{L^2(S^{d-1})}=-<\Delta_{S^{d-1}} Y_i^{(n)},Y_j^{(m)}>_{L^2(S^{d-1})}\\
&=-n(2-d-n)<Y_i^{(n)},Y_j^{(m)}>_{L^2(S^{d-1})}=n(n+d-2)\delta_{(i,n)(j,m)}\end{split}\ee
for $1\le n$ and $1\le i\le a_n$.
Therefore we have
\be\begin{split}
    \int_{S^{d-1}}|\frac 1r\nabla_{S^{d-1}}v|^2d\theta&=\int_{S^{d-1}}|\sum_{n=1}^\infty(A_nr^{n-1}+B_nr^{1-d-n})\sum_{i=1}^{a_n}B_i^{(n)}\nabla_{S^{d-1}}Y_i^{(n)}|^2\\
&=\sum_{n=1}^\infty n(n+d-2)(A_nr^{n-1}+B_nr^{1-d-n})^2\sum_{i=1}^{a_n}(B_i^{(n)})^2,
   \end{split}\ee
and
\be \begin{split}\int_{S^{d-1}}|\frac{\partial v}{\partial r}|^2d\theta&= \int_{S^{d-1}}|\sum_{n=1}^\infty(nA_nr^{n-1}+(2-d-n)B_nr^{1-d-n})\sum_{i=1}^{a_n}B_i^{(n)}Y_i^{(n)}|^2\\
&=\sum_{n=1}^\infty(nA_nr^{n-1}+(2-d-n)B_nr^{1-d-n})^2\sum_{i=1}^{a_n}(B_i^{(n)})^2.  \end{split}\ee
Consequently
\be \begin{split}&\int_{S^{d-1}}|\nabla v|^2d\theta=\\
&\sum_{n=1}^\infty[(nA_nr^{n-1}+(2-d-n)B_nr^{1-d-n})^2+n(n+d-2)(A_nr^{n-1}+B_nr^{1-d-n})^2]\sum_{i=1}^{a_n}(B_i^{(n)})^2.\end{split}\ee
Next
\be \begin{split}
&\int_{\overline B(0,1)\bs B(0,r_0)}|\nabla v|^2=\int_{r=r_0}^1 r^{d-1}dr\int_{S^{d-1}}|\nabla v|^2d\theta\\
=&\int_{r=r_0}^1[(nA_nr^{n-1}+(2-d-n)B_nr^{1-d-n})^2+n(n+d-2)(A_nr^{n-1}+B_nr^{1-d-n})^2]r^{d-1}dr\\
=&\int_{r=r_0}^1 [(2n^2+nd-2n)A_n^2r^{2n-2}+(2n+d-2)(n+d-2)B_n^2r^{2-2d-2n}]r^{d-1}dr\\
=&\int_{r=r_0}^1 (2n^2+nd-2n)A_n^2r^{2n+d-3}+(2n+d-2)(n+d-2)B_n^2r^{1-d-2n}dr\\
=&nA_n^2(1-r_0^{2n+d-2})+(n+d-2)B_n^2(r_0^{2-d-2n}-1).
   \end{split}\ee

Then by (7.17), we have
\be\begin{split} &nA_n^2(1-r_0^{2n+d-2})+(n+d-2)B_n^2(r_0^{2-d-2n}-1)\\
=&\frac{n(n+d-2)[r_0^{2-d-2n}-1][n+(n+d-2)r_0^{2n+d-2}]}{((n+d-2)r_0^n+nr_0^{2-d-n})^2}.\end{split}\ee
But $r_0<\frac12$, hence for all $n\ge 1$,
\be\frac{r_0^{2-d-n}-1}{r_0^{2-d-n}+1}\ge\frac{r_0^{-1}-1}{r_0^{-1}+1}\ge\frac{(\frac12)^{-1}-1}{(\frac12)^{-1}+1}=\frac13,\ee
therefore
\be\begin{split}
    &nA_n^2(1-r_0^{2n+d-2})+(n+d-2)B_n^2(r_0^{2-d-2n}-1)\\
\ge&\frac{n(n+d-2)[r_0^{2-d-2n}+1][n+(n+d-2)r_0^{2n+d-2}]}{3((n+d-2)r_0^n+nr_0^{2-d-n})^2}\\
\ge&\frac{n[nr_0^{2-d-2n}+n+d-2][n+(n+d-2)r_0^{2n+d-2}]}{3((n+d-2)r_0^n+nr_0^{2-d-n})((n+d-2)r_0^n+nr_0^{2-d-n})}\\
&=\frac n3[\frac{nr_0^{2-d-2n}+n+d-2}{(n+d-2)r_0^n+nr_0^{2-d-n}}][\frac{n+(n+d-2)r_0^{2n+d-2}}{(n+d-2)r_0^n+nr_0^{2-d-n}}]\\
=&\frac n3r_0^{d-2}.
   \end{split}
\ee

As a result,
\be\begin{split} \int_{B(0,1)\backslash B(0,r_0)}&|\nabla v|^2d\theta \ge\sum_{n=1}^\infty\frac n3r_0^{d-2}\sum_{i=1}^{a_n}(B_i^{(n)})^2\\
&\ge\frac13r_0^{d-2}\sum_{n=1}^\infty\sum_{i=1}^{a_n}(B_i^{(n)})^2=\frac13r_0^{d-2}||u_0(\theta)-m(u_0)||_{L^2(S^{d-1})}^2\\
&=\frac13r_0^{-1}\int_{\partial B(0,r_0)}|u_0-m(u_0)|^2.\end{split}\ee
Thus complete the proof of Proposition 7.8.\qed

\begin{cor}\label{cas1}Let $r_0>0,q\in\R^d$ be such that $r_0<\frac12 d(q,\partial B(0,1))$. Let $u_0\in C^1(\partial B(q,r_0)\cap \R^d,\R)$ and denote by $m(u_0)=(s_{d-1}r_0^{d-1})^{-1}\int_{\partial B(q,r_0)}u_0$ its average.

Then for all $u\in C^1((\overline{B(0,1)}\backslash B(q,r_0))\cap \R^2, \R)$ that satisfies
\be u|_{\partial B(q,r_0)}=u_0\ee
we have
\be\int_{B(0,1)\backslash B(q,r_0)}|\nabla u|^2\ge \frac13r_0^{-1}\int_{\partial B(q,r_0)}|u_0-m(u_0)|^2.\ee
\end{cor}

\nd The same as Corollary 7.23 of \cite{2p}.\qed

\begin{lem}
Let $d\ge 2$ be an interger, $\frac 14<r_0<1,$ $u\in C^1(B(0,1)\backslash B(0,r_0)\cap \R^d,\R)$ be such that $u|_{\partial B(0,r_0)}=\delta r_0,\ u|_{\partial B(0,1)}=0$; then
\be\int_{B(0,1)\backslash B(0,r_0)} |\nabla u|^2\ge c(d)\delta^2r_0^d,\ee
where $c(d)=\frac{(d-2)s_{d-1}}{\log 4}$. 
\end{lem}

\nd

For the case of $d=2$, please refer to \cite{2p} Section 7. Here we only prove the proposition for $d\ge 3$. 

Set $f(r,\theta)=A r^{-d+2}-A$ with $A=\frac{\delta r_0}{r_0^{-d+2}-1}$. Then $f$ is the harmonic function with the given boundary values. We have
\be\frac{\partial f}{\partial r}=(2-d)Ar^{1-d},\ \nabla_{S^{d-1}}f=0,\ee
hence
\be|\nabla f|^2=A^2(2-d)^2r^{2-2d}.\ee
Denote by $s_{d-1}=\H^{d-1}(S^{d-1})$, then
\be\begin{split}
\int_{B(0,1)\backslash B(0,r_0)} |\nabla f|^2&=\int_{S^{d-1}}d\theta\int_{r_0}^1r^{d-1}dr|\nabla f|^2\\
&=s_{d-1}\int_{r_0}^1r^{d-1}drA^2(2-d)^2r^{2-2d}\\
&=s_{d-1}A^2(2-d)^2\frac{1-r_0^{2-d}}{2-d}=\frac{s_{d-1}(d-2)\delta^2}{(r_0^{2-d}-1)}r_0^2\\
&\ge c(d)\frac{\delta^2r_0^2}{r_0^{2-d}-1}=\frac{r_0^{2-d}}{r_0^{2-d}-1}c(d)\delta^2r_0^d\ge c(d)\delta^2r_0^d.
\end{split}\ee
So we have
\be\int_{B(0,1)\backslash B(0,r_0)} |\nabla u|^2\ge c(d)\delta^2r_0^d\ee
since $f$ is harmonic.\qed

\begin{cor}\label{level}
For all $0<\epsilon<1$, there exists $C=C(\epsilon)>100$ such that if $0<r_0<1,$ $u\in C^1(\ B(0,1)\backslash B(0,r_0)\cap \R^d,\R)$ and
\be||u|_{\partial B(0,r_0)}-\delta||_\infty <\frac{\delta r_0}{C}\ et\ ||\ u|_{\partial B(0,1)}||_\infty <\frac{\delta r_0}{C},\ee
then 
\be\int_{B(0,1)\backslash B(0,r_0)} |\nabla u|^2\ge\epsilon c(d)\delta^2r_0^d.\ee
\end{cor}

\nd Apply Lemma 7.39 of \cite{2p}: take $r=r_0,\ f=u,$ and for each $C$, denote by $g$ the harmonic function with $g|_{\partial B(0,1)}=\frac\delta C r_0, g|_{\partial B(0,r)}=(1-\frac 1C) r_0\delta$. Then we have
\be\int_{B(0,1)\backslash B(0,r_0)} |\nabla u|^2\ge\int_{B(0,1)\backslash B(0,r_0)} |\nabla g|^2=(1-\frac 2C)^2c(d)\delta^2r_0^d.\ee
Thus for each $\epsilon<1$ we can always find a large enough $C$ such that $(1-\frac 2C)^2\ge\epsilon$.\qed

\section{Conclusion}

The rest of the proof of Theorem \ref{main} is almostly the same as the proof of the Almgren minimality of the union of two almost orthogonal planes in \cite{2p}. So we will only roughly describe what happens, without much detail. 

So fix $k$ large and $\e$ small enough. Set $D(x,r)=D_k(x,r),C^i(x,r)=C_k^i(x,r)$ for $1\le i\le m$, and $d_{x,r}=d_{x,r}^k$. For each $k$ fixed, we have chosen $o_k$ and $r_k$ as in the end of Section 5. Then by Proposition 6.1 (1), $E_k\cap D_k(0,\frac {39}{40})\bs D_k(o_k,\frac{1}{10}r_k)$ is composed of $m$ $C^1$ graphs $G^i,1\le i\le m$ on $P_k^1,\cdots, P_k^m$, such that (6.2) and (6.3) hold, where we replace $q_n,s_n$ by $o_k,r_k$. Moreover we can also suppose that $r_k<2^{-5}$, since $k$ is large. Inside $D_k(o_k,r_k)$, $E_k$ is $\e r_k$ far from any translation of $P_k$. By a compactness argument we know that the part $E_k\cap D(o_k,r_k)\bs D_k(o_k,\frac{1}{10}r_k)$ is composed of $m$ $C^1$ graphs, and is $\d r_k$  ($\d=\d(\e)$ comes naturally from the compactness argument, and only depends on $\e$) far from any translation of $P_k$ (cf. \cite{2p} Proposition 8.1 and Corollary 8.24). This makes one of the $m$ graphs $G^1,\cdots, G^m$, say $G^1$, being $\d r_k$ far from any translate of the $d-$plane $P_k^1$ in $D(o_k,r_k)\bs D_k(o_k,\frac{1}{10}r_k)$.

Then denote by $P=P_k^1$ for short, and let $g^1$ be as in (6.2); then $g^1$ is a map from $P$ to $P^\perp$, and is therefore from $\R^d$ to $\R^{(m-1)d}$. Write $g^1=(\varphi_1,\cdots, \varphi_{(m-1)d})$, where $\varphi_i:\R^d\to \R$, $1\le i\le (m-1)d$. Then since the graph of $g^1$ is $\d r_k$ far from all translation of $P$, there exists $1\le j\le (m-1)d$ such that
\be \sup_{x,y\in P\cap D(o_k, r_k)\backslash D(o_k, \frac 14 r_k)}|\varphi_j(x)-\varphi_j(y)|\ge C(m,d)r_k\delta,\ee
where $C(m,d)$ only depends on $m$ and $d$, for example we can take $C(m,d)=\frac{1}{\sqrt{(m-1)d}}.$

Suppose this is true for $j=1$. Denote by 
\be K=\{(z, \varphi_1(z)): z\in (D(0,\frac 34)\backslash D(o_k,\frac14 r_k))\cap P\},\ee 
then
\be \begin{split}K\mbox{ is the orthogonal }&\mbox{projection of }G^1\cap D(0,\frac 34)\\
&\mbox{ on a }d+1-\mbox{dimensional subspace of }\R^{md}.\end{split}\ee 


For $\frac 14r_k\le s\le r_k$, define 
\be \Gamma_s=K\cap p^{-1}(\partial D(o_k, s)\cap P)=\{(x,\varphi_1(x))|x\in\partial D(o_k, s)\cap P\}\ee  the graph of $\varphi_1$ on $\partial D(o_k, s)\cap P$. 

We know that the graph of $\varphi_1$ is $C(m,d)\d r_k$ far from $P$ in $D(o_k,r_k)\bs D(o_k,\frac 14r_k)$; then there are two cases:

1st case: there exists $t\in[\frac 14r_k, r_k]$ such that
\be \sup_{x,y\in \Gamma_t}\{|\varphi_1(x)-\varphi_1(y)|\}\ge \frac{\delta}{C}r_k,\ee 
where $C=4C(m,d)^{-1}C(\frac 12)$, $C(\frac 12)$ being the constant of Corollary \ref{level}.

Then there exists $a,b\in \Gamma_t$ such that $|\varphi_1(a)-\varphi_1(b)|>\frac\delta C r_k\ge\frac\delta Ct$. Since $||\nabla \varphi_1||_\infty\le ||\nabla \varphi||_\infty<1$, we have 
\be \int_{\Gamma_t}|\varphi_1-m(\varphi_1)|^2\ge t^{d+1}C'(\d),\ee
where $C'(\d)$ only depends on $d,m,\d$. 
Now in $D(0,\frac 34)$ we have $d(0,o_k)<6\e\le 10\e \cdot\frac34$, and $s<r_k<\frac 18<\frac 12\times\frac 34$, therefore we can apply Corollary \ref{cas1} and obtain
\be \int_{(D(0,\frac 34)\backslash D(o_k, t))\cap P}|\nabla \varphi_1|^2\ge C_1(\delta)t^d.\ee 

2nd case: for all $\frac 14r_k\le s\le r_k$, 
\be \sup_{x,y\in \Gamma_s}\{|\varphi_1(x)-\varphi_1(y)|\}\le \frac{\delta}{C}r_k.\ee 
However, since
\be\begin{split}
\frac12 r_k\delta&\le\sup\{|\varphi_1(x)-\varphi_1(y)|:x,y\in P\cap D(o_k,r_k)\backslash D(o_k,\frac14 r_k)\}\\
&=\sup\{|\varphi_1(x)-\varphi_1(y)|:s,s'\in[\frac14r_k,r_k],x\in\Gamma_s,y\in\Gamma_{s'}\},
\end{split}\ee
there existe $\frac 14r_k\le t<t'\le r_k$ such that
\be \sup_{x\in \Gamma_{t},y\in\Gamma_{t'}}\{|\varphi_1(x)-\varphi_1(y)|\}\ge C(m,d)r_k\delta. \ee 
Fix $t$ and $t'$, and without loss of generality, suppose that
\be \sup_{x\in \Gamma_{t},y\in\Gamma_{t'}}\{\varphi_1(x)-\varphi_1(y)\}\ge C(m,d)r_k\delta .\ee 
Then
\be \inf_{x\in\Gamma_{t}}\varphi_1(x)-\sup_{x\in\Gamma_{t'}}\varphi_1(x)\ge C(m,d)r_k\delta-2\frac\delta Cr_k=C(m,d)(1-\frac{1}{C(\frac12)})r_k\ge C(m,d)(1-\frac{1}{C(\frac12)})t'\ee 
because $C=4 C(\frac 12)$.

Now look at what happens in the ball $D(o_k,t')\cap P$. Apply Corollary \ref{level} to the scale $t'$, and since $\frac{t'}{t}\le 4, t'>t$, we get
\be \int_{(D(o_k,t')\backslash D(o_k, t))\cap P}|\nabla \varphi_1|^2\ge C(\delta,\frac 12,d,m)t'^d\ge C_2(\delta)t^d.\ee 

In both cases we pose $t_k=t$. The discussion above yields that there exists a constant $C_0(\delta)=\min\{C_1(\delta), C_2(\delta)\}$, which depends only on $\delta$, such that
\be \int_{(D(0,\frac 34)\backslash D(o_k, t_k))\cap P}|\nabla \varphi_1|^2\ge C_0(\delta)t_k^d.\ee 

On the other hand, since $|\nabla \varphi_1|\le|\nabla g^1|<1$,
\be \sqrt{1+|\nabla\varphi_1|^2}>\sqrt{1+\frac 12|\nabla\varphi_1|^2+\frac{1}{16}|\nabla\varphi_1|^4} =1+\frac14|\nabla\varphi_1|^2.\ee 
Hence 
\be\begin{split}
\H^d(K\bs C_k^1(o_k,t_k) )&=\int_{D(0,\frac 34)\backslash C_k^1(o_k,t_k)\cap P}\sqrt{1+|\nabla\varphi_1|^2}\ge\int_{D(0,\frac 34)\backslash C_k^1(o_k, t_k)\cap P}1+\frac14|\nabla\varphi_1|^2\\
&\ge \H^d((D(0,\frac 34)\backslash C_k^1(o_k, t_k))\cap P))+\frac 14\int_{D(0,\frac 34)\backslash C_k^1(o_k, t_k)\cap P}|\nabla \varphi_1|^2\\
&= \H^d((D(0,\frac 34)\backslash C_k^1(o_k, t_k))\cap P_k^1))+C_0(\d)t_k^d.
\end{split}\ee

Then by (8.3) we get
\be\begin{split} \H^d&(G^1\cap D(0,\frac34)\bs D(o_k,t_k))\ge \H^d(K\bs D(o_k,t_k))\\
&\ge \H^d((P_k^1+o_k)\cap D(0,\frac 34)\backslash D(o_k, t_k))+C(\d)t_k^2\\
&=\H^d(P_k^1\cap D(0,\frac 34)\backslash D(0, t_k))+C_0(\d)t_k^d,\end{split}\ee

 so that by the estimates in Section 7, we know that the measure of $G^1_{\frac{1}{10}r_k}$ is at least $C(\e)r_k^2$ more than its projection to $P_k^1$, that is, 
\be \H^d(G^1_{\frac{1}{10}r_k})\ge \H^d(D(0,1)\bs D(o_k,\frac{1}{10}r_k)\cap P_k^1)+C(\e)r_k^d.\ee

For the other $2\le i\le m$, $G^i$ is a graph on $D(0,1)\bs D(o_k,\frac{1}{10}r_k)\cap P_k^i$, thus 
\be \H^d(G^i_{\frac{1}{10}r_k})\ge \H^d(D(o_k,\frac{1}{10}r_k)\cap P_k^i), 2\le i\le m.\ee

For the part of $E_k$ inside $D(o_k,\frac{1}{10}r_k)$, by hypotheses the characteristic angle $\a(k)$ of $P_k$ is larger than $\frac\pi 2-\frac 1k$, hence Lemmas 2.16 and 2.27 gives
\be \H^d(E_k\cap D(o_k,\frac{1}{10}r_k))\ge (1+2\cos(\frac\pi2-\frac1k))^{-1}\H^d(P_k\cap D(o_k,\frac{1}{10}r_k)).\ee

We sum over (8.18)-(8.20), and get
\be \H^d(E_k)\ge \H^d(P_k\cap D(0,1))+[C(\e)-C\cos(\frac\pi2-\frac1k)]r_k^d,\ee
where $C(\e)$ depends only on $\e$. Thus for $k$ large enough, (8.21) gives
\be \H^d(E_k)>\H^d(P_k\cap D(0,1))=mv(d),\ee
 which contradicts (4.10). So we get the desired contradiction.
 
 \section{Corollary and some open questions}
 
 \begin{defn} We call that $\a=(\a^ij)_{1\le i<j\le m}$, where $\a^{ij}=(\a^{ij}_1,\cdots, \a^{ij}_d)$ ($\a^{ij}_1\le\cdots\le\a^{ij}_d$) is a topological (resp. Almgren) minimal angle if the union of $m$ $d-$dimensional planes with characteristic angle $\a$ is topologically (resp. Almgren) minimal.
 \end{defn}
 
\begin{cor}Let $d\ge 3$, then for each $k\le d-2$ and $m\in \N$, there exists $\theta=\theta(d,k,m)\in]0,\frac\pi 2[$, such that if $\a=(\a^ij)_{1\le i<j\le m}$ satisfies $\a_l^{ij}=0$ for $1\le l\le k,1\le i,j\le m$, $\a_l>\theta(d,k,m)$ for $l>k,1\le i,j\le m$, then $\a$ is a topological minimal angel.\end{cor}

\nd By Theorem \ref{main}, take $\theta=\theta_{m,d-k}$, then for any $m$ $d-$planes as in the theorem, the angle condition 
\be \a_l^{ij}=0 \mbox{ for }1\le l\le k,1\le i,j\le m, \a_l>\theta(d,k,m)\mbox{ for }l>k,1\le i,j\le m\ee
means that $R=\cap_{i=1}^m P^i$ is a $k$-dimensional plane, and $P^i=R\times Q^i$, where $Q_i\perp R$ is a $d-k$dimensional plane, $1\le i\le m$. The characteristic angles between $Q^i$ and $Q^j$ are just $\a_{k+1}^{ij},\cdots \a^{ij}_d$, which are larger than $\theta(m,d-k)$, hence their union $\cup_{i=1}^m Q^i$ is topologically minimal. As a result, by Proposition 3.23 of \cite{topo}, $\cup_{i=1}^m P^i=R\times \cup_{i=1}^m Q^i$ is also topologically minimal. \qed

\begin{rem} Notice that in Corollary 9.1, for difference $k$, the families of unions of planes are difference, moreover they are far from each other. 

It is interesting to ask whether we have some kind of "interpolation" type property, that is, if $\a=(\a^{ij}_k)_{1\le k\le d,1\le i,j\le m}$ and $\beta=(\beta^{ij}_k)_{1\le k\le d,1\le i,j\le m}$ are minimal angles, and $\gamma=(\gamma^{ij}_k)_{1\le k\le d,1\le i,j\le m}$ satisfies $\a^{ij}_k\le \gamma^{ij}_k\le \beta^{ij}_k$ for all $i,j,k$, then is $\gamma$ a minimal angle? However, we even can not answer the following simpler question:

If $\a=(\a^{ij}_k)_{1\le k\le d,1\le i,j\le m}$ is a minimal angle, and $\beta=(\beta^{ij}_k)_{1\le k\le d,1\le i,j\le m}$ is larger than $\a$ (i.e. $\a^{ij}_k\le \beta^{ij}_k$ for all $i,j,k$), is $\beta$ a minimal angle? 

Intuitively, the bigger the angles are, the more likely is the union of two planes to be minimal. But the above question is still open.\end{rem}

\renewcommand\refname{References}
\bibliographystyle{plain}
\bibliography{reference}

\end{document}